\documentclass[11pt]{amsart}
\usepackage{amsmath}
\usepackage{amssymb}
\usepackage{euscript}
\usepackage{pifont}
\usepackage{tikz}
\usepackage[all]{xy}
\numberwithin{equation}{section}
\numberwithin{figure}{section}
\input cyracc.def

\usepackage{amssymb,amsmath}
\usepackage{pstricks}
\usepackage{eepic}

\newrgbcolor{blue}{0 0 .7}

\newtheorem{definition}[equation]{Definition}
\newtheorem{theorem}[equation]{Theorem}

\newtheorem{lemma}[equation]{Lemma}
\newtheorem{example}[equation]{Example}
\newtheorem{remark}[equation]{Remark}

\newcommand{\w}{\wedge}
\newcommand{\n}{\notag}
\newcommand{\noi}{\noindent}
\newcommand{\im}{\mbox{i}}
\newcommand{\hook}{\hookrightarrow}

\newcommand{\vsp}{\vspace}
\newcommand{\sq}{ \;  \square}
\newcommand{\np}{\newpage}
\newcommand{\lag}{\langle}
\newcommand{\rag}{\rangle}

\newcommand{\E}{ \mathbb{E}}
\newcommand{\C}{ \mathbb{C}}
\newcommand{\R}{ \mathbb{R}}

\newcommand{\bomega}{\overline{\omega}}

\newcommand{\bh}{\overline{h}}

\newcommand{\bpartial}{\overline{\partial}}

\newcommand{\x}{\textnormal{x}}

\newcommand{\ES}{\textnormal{S}}
\newcommand{\SO}{\textnormal{SO}}

\newcommand{\JAI}{\textnormal{J}}

\newcommand{\Pf} {\noi \emph{Proof.}$\; \;$}

\newcommand{\beq}{ \begin{equation}}
\newcommand{\eeq}{ \end{equation} }

\makeatletter

\newcommand{\Rmnum}[1]{\expandafter\@slowromancap\romannumeral #1@}
\makeatother

\begin{document}

\sloppy

\title{Degree 3 algebraic minimal surfaces in the 3-sphere}
\author{Joe S. Wang}
\email{jswang12@gmail.com}
\begin{abstract}
We give a local  analytic characterization
that a minimal surface in the 3-sphere $\, \ES^3 \subset \R^4$
defined by an irreducible cubic polynomial
is one of the Lawson's minimal tori.
This provides  an alternative proof of the result  by  Perdomo
(\emph{Characterization of order 3 algebraic immersed minimal surfaces of $S^3$},
Geom. Dedicata  129  (2007),   23--34).

\end{abstract}
\maketitle
\tableofcontents

\section{Introduction}\label{sec1}
Let $\, \E=\R^4$ be the four dimensional Euclidean vector space.
Let $\, \ES^3 \subset \E$ be the unit sphere equipped with the induced Riemannian metric.
A \emph{minimal surface}  $\, M \hook \ES^3$ is by definition
an immersed surface with vanishing mean curvature,
which locally gives an extremal for the area functional.
The differential equation for minimal surfaces
is an elliptic Monge-Ampere  equation
defined on the unit tangent bundle of $\, \ES^3$.
It is well known that
the  minimal surface  equation   is  locally equivalent to
the  elliptic  $\, \sinh$-Gordon equation
for  one  scalar function of two variables;
\beq\label{1sinh}
\Delta u + \sinh{u} = 0.
\eeq

Given a  surface $\, M \hook \ES^3$,
consider the cone $\, O   \ast     M \hook \E$ over the origin $\,  O \in \E$.
Then $\, M \hook \ES^3$ is minimal
whenever
$\, O  \ast    M \hook \E$ is  minimal   as a hypersurface in  the Euclidean space $\, \E$.
In this respect,
a minimal surface $\, M \hook \ES^3$ is called \emph{algebraic of degree  m}
if  there exists a nonzero  irreducible degree $\, m$   homogeneous polynomial $\, F: \E \to \R$
which vanishes on $\, M$.
A nonzero   irreducible   homogeneous polynomial  $\, F$
defines a (possibly singular) minimal surface   when
\beq\label{1poly}
 | \nabla F|^2 \, \Delta F  -  \, (\nabla F)^t \cdot  H(F) \, \cdot  (\nabla F) \equiv 0, \; \mod \; F.
\eeq
Here $\, \nabla F,\, \Delta F,  \, H(F)$  denote
the gradient, the Laplacian, and the Hessian of the function $\, F$ respectively
with respect to the Euclidean metric of $\, \E$, \cite[p260]{Hs}.

Hsiang  classified the    homogeneous minimal hypersurfaces
in the standard  Euclidean spheres, \cite{Hs}.
In \cite{Hs},
he also  showed  that
the totally geodesic 2-sphere  and Clifford torus
are the only algebraic minimal surfaces  in $\, \ES^3$ of degree $\, \leq   2$.
Lawson  constructed  an infinite sequence of
algebraic minimal tori(or Klein bottles)  in $\, \ES^3$ of arbitrary high degree, \cite{La}.
Hsiang and Lawson  gave an analysis for
the  minimal submanifolds of low cohomogeneity in general Riemannian homogeneous spaces,
where, in particular,
the Lawson's sequence of algebraic minimal tori
are extended to a countable family of  cohomogeneity 1 minimal tori in $\, \ES^3$, \cite[p32]{HsL}.

Recently,
Perdomo gave a  characterization  of degree 3 algebraic minimal surfaces in $\, \ES^3$
as  one of  the Lawson's algebraic minimal tori,  \cite{Pe05, Pe}.
One of the main idea of  his analysis    is
that such a minimal surface  necessarily contains   a great circle,
which puts the defining cubic polynomial into a special normal  form.
The minimal surface equation  \eqref{1poly} is then applied
to  successively normalize the polynomial coefficients by  differential  algebraic analysis.

The purpose of the present paper
is to give an  analytic characterization of  degree 3 algebraic minimal surfaces in $\, \ES^3$.
We employ the method of local differential analysis
and show that
Lawson's  algebraic minimal torus of degree 3
is the only minimal surface in $\, \ES^3$
which satisfies the compatibility equation
to lie in the  zero locus of an irreducible  cubic polynomial.

\textbf{Main results.}

1.
The structure equation
for the degree 3 algebraic minimal surfaces in $\,\ES^3$ is determined, Theorem \ref{4thm}.
The analysis for the structure equation
provides  the formula for the defining cubic polynomial,
and it  is explicitly identified  as  the Lawson's degree 3 example.

2.
The structure equation shows that
the degree 3 algebraic minimal surface
is  the  conjugate surface of   a  principally bi-planar\footnotemark
\, minimal surface.
It has the Killing nullity 1,
and the curvature takes values in the closed interval $\, [\, -3, \, \frac{3}{4}\, ]$.

\footnotetext{
A minimal surface in $\, \ES^3$ is \emph{principally bi-planar}
when each of  its principal curves lies in a  totally geodesic  $\, \ES^2 \subset  \ES^3$, \cite{Yam}.
There exist  locally one  parameter family of  distinct
principally bi-planar minimal surfaces in $\,  \ES^3$.}

Let us give an outline of the analysis.
The problem of  finding  an algebraic  minimal surface in $\, \ES^3$
can be interpreted as
the problem of finding a solution to  the $\, \sinh$-Gordon equation  \eqref{1sinh}
with the property that
it  admits an associated, in a certain algebraic manner,  polynomial $\, F$
which  satisfies \eqref{1poly}.
This  is equivalent to finding
a constant section of an appropriate   tensor bundle over the minimal surface
which    vanishes  when evaluated on the surface.
As the degree of $\, F$ increases,
this  problem  imposes  a sequence of compatibility equations
on the higher order  structure functions of the minimal surface.
For the degree 3 case treated in this paper,
the compatibility equations are reduced to a pair of  third order equations.
\textbf{Main results} are obtained by
applying the  over-determined PDE analysis to  these equations.

The present work can be considered as
a coordinate free and equivariant interpretation of
the aforementioned Perdomo's original characterization.
Compared to his   analysis,
one possible advantage  would be that
this interpretation of   an algebraic  minimal surface
is susceptible to  local analytic method,
and could be, in theory,   applied  to general  higher degree cases.
However,
it  is  practically difficult to perform
the required  computations manually.
Our  analysis  relies essentially  on the use of  computer  machine.
On the other hand,
our analysis also suggests as a byproduct
a distinguished class of minimal surfaces
defined by an additional fourth order equation
($\, \Delta_4=0$ in \eqref{41Delta}).
An analysis of this class of  minimal surfaces will be reported
in the subsequent paper.

We carried out the analysis for the degree 4  algebraic minimal surfaces
by the same method used in this paper.\footnotemark
\,  The  analysis    indicates that
such a minimal surface   is
one of the family of  cohomogeneity 1 minimal tori described by Hsiang and Lawson.
Partly due to the length and complexity of the algebraic analysis involved,
we do not   present the details of  the  analysis for degree 4 case  in this paper.

\footnotetext{
The analysis of the degree 4  algebraic minimal surfaces
was the original motivation for the present work.
One needs a characterization of  the degree 3  surfaces  first
to exclude them from the analysis for the degree 4  surfaces.}

In a different context,
similar analysis can be applied to
special Legendrian surfaces in the 5-sphere,
which are the links of special Lagrangian 3-folds in $\, \C^3$, \cite{Ha}.
A special Legendrian surface is called
\emph{quasi-algebraic of degree $\,m$}
if it lies in the zero locus of
a nonzero  irreducible degree $\, m$ real homogeneous polynomial
$\, F: \C^3 \to \R$.
The analysis shows that;

1.
 a quasi-algebraic degree 2 special Legendrian surface
is necessarily
one of the  cohomogeneity 1 surfaces treated by Haskins, \cite{Ha},

2.
 a quasi-algebraic degree 3 special Legendrian surface
with an appropriate $\, \mathbb{Z}_3$ symmetry
is necessarily
a quasi-algebraic surface of degree 2.

The computation indicates  in fact   that
any quasi-algebraic degree 3 special Legendrian surface
is necessarily
a quasi-algebraic surface of degree 2.
We do not have a complete proof of this claim at this time of writing.

The paper is organized as follows.
In Section \ref{sec2},
we set  the basic structure equation  for the minimal surfaces in $\, \ES^3$,
and  compute its infinite sequence of prolongations, \eqref{21ecstruct}, Lemma \ref{211lemma}.
In order to simplify the computations,
the induced complex structure on the minimal surface is utilized,
and  the structure equations are  written in complex form.
In Section \ref{sec3},
we give a description of the method
to detect the compatibility equations for a minimal surface to be algebraic.
The method is applied to two simple cases of degree 1, and degree 2  algebraic minimal surfaces,
Theorem \ref{32thm},  Theorem \ref{33thm}.
In Section \ref{sec4},
we present the results of  differential analysis for  degree 3 algebraic  minimal surfaces,
Theorem \ref{4thm}.
After a preliminary reduction in Section \ref{sec41},
the analysis is divided into two cases.
The differential analysis shows that
a degree 3 algebraic minimal surface necessarily satisfies
an auxiliary equation either of  order 3, Section \ref{sec42},
or  of order  4, Section \ref{sec44}.
The auxiliary equation  of  order 3 is compatible
and supports a defining  irreducible cubic polynomial,
which is explicitly verified as the one  given by Lawson.


The majority of computations were performed using the computer algebra system
\texttt{Maple} with  \texttt{difforms}  package.\footnotemark

\footnotetext{
The  \texttt{Maple} worksheet is available upon request by email.}

Throughout the paper,
a surface is a connected  smooth two dimensional manifold.

\section{Minimal surfaces in $\, \ES^3$}\label{sec2}
In this section,
we  set   the  basic structure equations  for a minimal surface  in  $\, \ES^3$.
In Section \ref{sec21},
we apply the  moving frame method
to  determine the structure equation
for   an immersed  oriented  minimal surface  in $\, \ES^3$.
In Section \ref{sec211},
we compute the infinite prolongation  of the structure equation explicitly,
Lemma \ref{211lemma}.

The analysis  and  results in this section are classical, and well known.
For the standard reference on the  theory of minimal surfaces in  $\, \ES^3$,
we refer to \cite{La}  and the references therein.

The basic structure equations  established in this section
will be used  implicitly throughout the paper.

\subsection{Structure equation}\label{sec21}
Let $\, \E = \R^{4}$ be the four dimensional Euclidean vector space.
Let $\, \ES^3 \subset \E$ be the  unit sphere equipped with the induced Riemannian metric.
The special orthogonal group $\, \SO_4$ acts transitively on $\, \ES^3$
as a  group of  isometry,  and $\, \ES^3 = \SO_4 / \SO_3$.

Let $\, \Lambda \to \ES^3$ be the $\, \ES^2$-bundle of unit tangent vectors.
Let $\, Gr^+(2, \, \E)$ be the Grassmannian   of  two dimensional  oriented   subspaces of $\, \E$.
$\SO_4$ acts transitively on both $\, \Lambda$ and $\, Gr^+(2, \, \E)$,
and there exists the incidence double fibration;

\begin{picture}(300,109)(-67,-21) \label{21double}
\put(157,72){$ \SO_4$}
\put(172,57){$\pi$}
\put(162,57){$\downarrow$}
\put(161,40){$\Lambda = \SO_4 / \SO_2$}
\put(175,25){$\searrow$}
\put(145,25){$\swarrow$}
\put(184,10){ $Gr^+(2, \, \E)= \SO_4 / (\SO_2 \times \SO_2)$ }
\put(129,10){ $ \ES^3$ }
\put(187,27){$\pi_1$}
\put(133,27){$\pi_0$}
\put(99, -11){ \textnormal{Figure \ref{21double}.  Double fibration }}
\end{picture}

To fix the notation once and for all,
let us define the projection maps $\, \pi,  \, \pi_0$, and  $\, \pi_1 \, $ explicitly.
Let  $ e= (e_0, \, e_1, \, e_2, \, e_3)$
denote   the $\, \SO_4 \subset $ GL$_4\R$\, frame of $\, \E$.
Define
\begin{align}\label{21defipi}
\pi (e) &= ( e_0, \,  e_0 \w e_3),   \\
\pi_0 ( e_0, \,  e_0 \w e_3) &=e_0,  \n \\
\pi_1 ( e_0, \,  e_0 \w e_3) &= e_0 \w e_3. \n
\end{align}

The $\, \SO_4$-frame $\, e$  satisfies the structure equation
\begin{align}\label{21frame}
d e_A &= \sum_B \,e_B \, \omega^B_A,   \\
\omega^A_B+\omega^B_A &= 0, \n
\end{align}
for the Maurer-Cartan form $\, (\, \omega^A_B)$  of  $\, \SO_4$.
$\, (\, \omega^A_B)$ satisfies the compatibility   equation
\beq\label{2MC}
d\omega^A_B+ \sum_C \omega^A_C \w \omega^C_B=0.
\eeq

Let $\, \x : M \hook \ES^3$ be an   immersed     oriented surface.
By the general  theory of moving frames,
there exists a  lift
$\, \tilde{\x}: M \hook \Lambda$
such that $\, e_0 = \x$,  and $\, e_3$ is the oriented  normal to the surface $\, \x$.
Let $\, \tilde{\x}^* \SO_4 \to M$ be the pulled back $\, \SO_2$-bundle.
We continue to use $\, (\omega^A_B)$ to denote the  pulled back
Maurer-Cartan form  on $\, \tilde{\x}^* \SO_4$.
From  \eqref{21defipi}, \eqref{21frame},
the initial state of   $\, (\omega^A_B)$ on  $\, \tilde{\x}^* \SO_4$ takes the form
\beq\label{21omega}
(\omega^A_B) =
\begin{pmatrix}
\cdot        & -\omega^1 & -\omega^2 & \cdot \\
\omega^1 &  \cdot        &  \omega^1_2  & \omega^1_3 \\
\omega^2 &  \omega^2_1 & \cdot         & \omega^2_3 \\
\cdot        &  \omega^3_1 & \omega^3_2 & \cdot
\end{pmatrix}.
\eeq
Here  '$\cdot$' denotes zero, and  $\, \omega^A_0 = \omega^A,  \, A = 1, \, 2$.
By definition,
$\, \textnormal{\Rmnum{1}} = \lag \, de_0, \,  de_0 \, \rag
              = (\omega^1)^2 + (\omega^2)^2$
is the induced Riemannian metric of the immersed surface,
where $\, \lag \,  \,, \, \rag$ is the inner product of $\, \E$.

Differentiating $\, \omega^3_0=0$,
one gets
\beq
\omega^3_1 \w \omega^1 + \omega^3_2 \w \omega^2 = 0. \n
\eeq
By Cartan's lemma,
there exist coefficients $\,  h_{AB},   \, A, B =1,2$,  symmetric in   indices such that
\beq
\omega^3_A  = \sum_B h_{AB} \, \omega^B. \n
\eeq
The structure equation shows that
the quadratic differential
\beq\label{21defiII}
\textnormal{\Rmnum{2}} = \omega^3_A \circ \omega^A
                       =  h_{AB}  \, \omega^A  \circ  \omega^B
\eeq
is well defined on $\, M$.
\Rmnum{2} is   the second fundamental form of  the immersed surface $\, \x$.
\begin{definition}
Let $\, \x: M \hook \ES^3$ be an immersed  oriented surface.
Let  $\, \textnormal{II}$ be the quadratic differential  \eqref{21defiII}
which is  the second fundamental form of the immersed surface $\, \x$.
$\, \x$   is  \emph{minimal}
when the trace of  the quadratic  differential  $\, \textnormal{\Rmnum{2}}$
with respect to the induced metric $\, \textnormal{\Rmnum{1}}$ vanishes,
or equivalently when
\beq
h_{11}+h_{22}=0. \n
\eeq
\end{definition}
\noi
From now on,
a minimal surface would mean an immersed  oriented   minimal surface.

In order to utilize the induced complex  structure on $\, M$ as a Riemann surface,
let us  introduce the  complexified structure equation.
Let $\, \E^{\C} = \E \otimes \C$ be  the complexification,
and
consider  the following  $\, \E^{\C}$-frame.
\begin{align}\label{21ecframe}
e^{\C} &=(e_0, \, E_1, \, E_{-1}, \, e_3 ),   \\
E_1 &= \frac{1}{2} ( e_1 - \im \, e_2), \quad \; \; \im^2=-1, \n \\
E_{-1}&=\overline{E}_1. \n
\end{align}
Rewriting \eqref{21frame}, \eqref{21omega} with respect to $\, e^{\C}$,
one gets
\beq
d  e^{\C} =  e^{\C}
\begin{pmatrix}\label{21ecstruct}
\cdot        & -\frac{1}{2}\bomega  &  -\frac{1}{2}\omega  & \cdot \\
\omega  &  -\im \, \rho      &   \cdot   &  - \bh_2 \bomega \\
\bomega &   \cdot       &  \im \, \rho         &  -  h_2 \, \omega  \\
\cdot        &  \frac{1}{2} h_2 \, \omega & \frac{1}{2} \bh_2 \bomega& \cdot
\end{pmatrix},
\eeq
where we set
\begin{align}
\omega & = \omega^1 + \im \, \omega^2, \n \\
\rho    &= \omega^1_2, \n  \\
h_2 & = h_{11} - \im \, h_{12}. \n
\end{align}
Differentiating \eqref{21ecstruct}, one gets the compatibility equations
\begin{align}\label{21struct}
d \omega&= \im \, \rho \w \omega,  \\
d \rho &= K \frac{\im}{2} \, \omega \w \bomega,
\; \; \mbox{where} \;\;  K= 1 - h_2 \bh_2, \n \\
dh_2 + 2   \, \im  \, h_2 \rho &\equiv 0, \mod\;  \omega. \n
\end{align}
Here
$\, K$ is the curvature of the induced metric on the minimal surface.
\begin{remark}
The subscript '2' in the notation '$\,h_2$'
represents the weight of the action  by the structure group $\, \SO_2$.
It  is convenient for a computational purpose.
\end{remark}

\eqref{21struct} shows that
the quadratic differential
\beq\label{21defiIIC}
\textnormal{\Rmnum{2}}^{\C}  = h_2 \, \omega \circ \omega
\eeq
is well defined on $\, M$, and is holomorphic
with respect to the complex structure on $\, M$ defined
by the $\, (1,  0)$-form $\, \omega$.
$\, \textnormal{\Rmnum{2}}^{\C}$ is
the complexified second fundamental form  of  the minimal surface.

\begin{definition}
Let $\, \x: M \hook \ES^3$ be an immersed oriented minimal surface.
The holomorphic quadratic differential
$\, \textnormal{\Rmnum{2}}^{\C}$, \eqref{21defiIIC},
is the \emph{Hopf differential} of the minimal surface.
The zero set of  $\, \textnormal{\Rmnum{2}}^{\C}$ is the \emph{umbilic divisor}.
\end{definition}
\begin{example}\label{21sphere}
Consider   an immersed minimal sphere   $\,   \ES^2 \hook  \ES^3$.
Since $\, \ES^2$ supports no nonzero holomorphic differentials of positive degree,
$\, \textnormal{\Rmnum{2}}^{\C} =0$
and the structure coefficient $\, h_2$ vanishes identically.
A minimal sphere in $\, \ES^3$ is   necessarily totally geodesic.

Conversely, it is clear from the structure equation \eqref{21ecstruct}
that
a minimal surface in $\, \ES^3$ with $\, h_2 \equiv 0$
is locally equivalent to the  totally geodesic sphere.
\end{example}

For   a compact   minimal  surface   $\, M \hook  \ES^3$
of  genus $\, g \geq 1$,
the umbilic divisor has   degree  $\, 4g-4$ by Riemann-Roch theorem.
In particular, the umbilic divisor of a minimal torus is empty.

\subsection{Prolongation}\label{sec211}

In this section,
we  compute the infinite sequence of  prolongations of the structure equation \eqref{21struct}.
The prolonged structure equation  will be   used implicitly
for  the differential analysis in  Section \ref{sec4}.

Let us introduce the higher order derivatives
of  the structure coefficient $\, h_2$  in \eqref{21struct} inductively by
\beq
dh_j +  \im   \, j    h_j  \, \rho
= h_{j+1} \, \omega + h_{j, -1} \, \bomega, \; \; j = 2, \, 3, \, ... \, .\n
\eeq
For a notational purpose,
denote
$\, h_{-j}= \bh_j, \; j = 2, \, 3, \, ... \,$.
For instance,
the curvature of the surface is written as
$\, K = 1 - h_2   \bh_2 = 1-h_2  h_{-2}$.

\begin{lemma}\label{211lemma}
\begin{align}
h_{2,-1}&=0, \,   \n \\
h_{j+1, -1}&=   \sum_{s=0}^{j-2}
c_{js} \,  h_{j-s} \, \partial^{s} K,
\; \; \; \mbox{for} \; \; j \geq 2,  \; \; \mbox{where} \n \\
& \quad  c_{js} =\frac{(j+s+2) }{2} \frac{ (j-1) !}{(j-s-2)!(s+2)!}
             =\frac{(j+s+2) }{2 j} {j \choose s+2}. \n
\end{align}
Here by definition
$\, \partial^{s} K =\delta_{0s} -h_{2+s} h_{-2}\,$  for $\, s \geq 0$.
\end{lemma}
\Pf
Differentiating $\, d h_j$ and collecting $\, \omega \w \bomega$-terms,
one gets
\beq
h_{j+1,-1} = \partial h_{j,-1} + \frac{j}{2}\, h_j \, K,  \n
\eeq
where  $\, \partial h_{j,-1}$ is the $\, \omega$-component of $\, d h_{j,-1}$
($\, \partial^s K$  is    defined similarly).
The formula is verified by direct computation.
Note that
\begin{align}
a_{j0}&=\frac{(j+2)(j-1)}{4}, \n \\
a_{j(j-2)}&=1.  \n
\quad \sq
\end{align}

 \begin{example}\label{211torus}
Consider   $\, \R^4 = \C \oplus \C = \C^2$.
Clifford torus is  the  minimal surface  in $\, \ES^3$
defined as the product of  circles
\beq
\{ \, (z_1, \, z_2)  \in  \C^2 \; | \quad    |z_1|^2 = |z_2|^2 = \frac{1}{2} \; \}. \n
\eeq
The induced Riemannian metric is flat,
and the curvature  $\, K$ vanishes.
Lemma \ref{211lemma} shows that
\beq
0 = \partial K = - h_3 h_{-2}. \n
\eeq
Clifford torus is not totally geodesic,
and this implies $\, h_3 \equiv 0$.

Conversely, consider a minimal surface in $\, \ES^3$
such  that the structure function  $\, h_3 \equiv 0$.
By Lemma \ref{211lemma},
$\, \bpartial h_3 = h_2 K  \equiv 0$,
and either $\, h_2 \equiv 0$(and $\, K=1$), or $\, K=0$.
In the latter case,
it is well known that the surface is locally congruent to Clifford torus.
\end{example}

\section{Algebraic minimal surfaces}\label{sec3}
Let $\, S^m(\E)$  be
the vector space of  real,  homogeneous degree $\,m$ polynomials on
the four dimensional Euclidean vector  space $\, \E$.
By the metric duality, we identify
$\, S^m(\E)$
with the space of symmetric $\, m$-tensors $\, Sym^m(\E)$.
Let $\, S(\E) =\oplus_{m=0}^{\infty} \,S^m(\E)$.
\begin{definition}\label{3defalg}
Let $\, \x: M \hook \ES^3 \subset \E$ be a   minimal surface in the unit sphere.
$\, \x$ is \emph{algebraic}
if there exists a nonzero   homogeneous polynomial  $\, F \in S(\E) $
which  vanishes on $\, \x$.
For each  $\, m \geq 0$,
let $\, \mathcal{J}^m_{\x} \subset S^m(\E)$
be the subspace  of   degree $\, m$ polynomials vanishing on $\, \x$.
\emph{Degree} of  an algebraic minimal surface $\, \x$ is
the minimum integer $\, m$ such that $\, \mathcal{J}^m_{\x}$ is nontrivial.
\end{definition}
\noi
Note by definition that for  an algebraic minimal surface  $\, \x$  of  degree $\, m_0$,
$\, \mathcal{J}^m_{\x}$ is nontrivial for all $\, m \geq m_0$.

In this section,
we first  describe the    idea of
how the over-determined  PDE analysis can be applied to
detect the compatibility equations for a minimal surface to be algebraic, Section \ref{sec31}.
We then apply this   to
the cases of degree 1, and degree 2  algebraic minimal surfaces,
Section \ref{sec32}, Section \ref{sec33}.

We continue to use the structure equations established in Section \ref{sec2}.
The method of  local  differential analysis  described  in this section
will be applied to degree 3  algebraic minimal surfaces in Section \ref{sec4}.

\subsection{Structure equation}\label{sec31}
Let $\, \x: M \hook \ES^3$ be a   minimal surface
with the associated    complexified  frame $\, e^{\C}$,  \eqref{21ecframe}.
Consider the trivial bundle $\, S^m(\E) \times \ES^3 \to \ES^3$,
and the induced bundle
$\, \x^*(S^m(\E) \times \ES^3) \to M$.
A section of $\, \x^*(S^m(\E) \times \ES^3)$ over $\,M$
can be represented
in terms of the frame $\, e^{\C}$ as follows.
\beq
F = \sum^{i+j+k+l=m}_{i, \, j, \, k, \, l \, \geq 0} p_{i j k l} e_0^i e_3^j E_1^k E_{-1}^l, \n
\quad p_{i j k l} = \overline{  p_{i j l k} },
\eeq
for a set of coefficients $\, p_{i j k l}$.
Equivalently,
such $\, F$ is an $\, S^m(\E)$-valued function on $\, M$.

The condition that
$\, F$ is a constant section, or $\, F \in  S^m(\E)$ with an abuse of notation,
is expressed by
\beq\label{31dF}
d F = 0.
\eeq
Expanding  the exterior derivative $\, d F$  by Leibniz rule
in the variables $\, \{ \, e_0, \, e_3, \, E_1, \, E_{-1} \, \}$,
one gets  the structure equation for the coefficients $\, \{ \, p_{i j k l} \, \}$;
\beq\label{31dpstruct}
d p_{i j k l} = \sum \phi_{i j k l}^{i' j' k' l'} p_{i' j' k' l'},
\eeq
where  $\, \phi_{i j k l}^{i' j' k' l'}$ is determined as
a linear combination of the components of  the Maurer-Cartan form  \eqref{21ecstruct}.

\begin{remark}
The  formulae for $\,\{ \, \phi_{i j k l}^{i' j' k' l'} \, \}$ will be given explicitly
in Section \ref{sec32}, Section  \ref{sec33}, Section  \ref{sec4}
for the cases $\, m=1, \, 2, \, 3$ respectively.
\end{remark}

The condition that the polynomial $\, F$ vanishes  on $\, \x$ is expressed by
\beq\label{31initial}
p_{m000}=0.
\eeq
The idea is then to apply the structure equation  \eqref{31dpstruct}  repeatedly
starting from  the initial state \eqref{31initial}
by using  \eqref{21struct} and \eqref{211lemma}.
As one differentiates,
the successive  derivatives of  \eqref{31initial} imply
more and more linear compatibility equations on the polynomial coefficients $\,\{ \,  p_{i j k l} \, \}$,
thereby reducing the number of independent coefficients.

In the case $\, F$ is irreducible,
the reduction process eventually leads to a single independent polynomial coefficient,
and  the  integrability equations  for a minimal surface to  support  such an irreducible  polynomial
are expressed  as a set of algebraic equations on the structure  functions
$\, \{ \, h_{\pm2}, \, h_{\pm 3}, \, ... \, \}$.
For the  minimal surface  satisfying these integrability equations,
the formula for   $\, F$ is  obtained  by
evaluating the polynomial coefficients
at an appropriate  generic point of the minimal surface.

In  the next two sections,
we examine the well known cases  of  degree 1, and  degree 2 minimal surfaces
following this method just described.
These two cases not only serve as exercise
for  the differential analysis  of the more complicated case
of   degree 3  surfaces in Section \ref{sec4},
but they are also  necessary prerequisites.
From the remark below Definition \ref{3defalg},
one needs the characterization of  degree 1, and degree 2 algebraic minimal surfaces
to exclude them from the  analysis for
the degree 3 surfaces.

\subsection{Degree  1  algebraic minimal surfaces}\label{sec32}
In this section,
we give a  characterization of degree 1 algebraic minimal surfaces.
\begin{theorem}\label{32thm}
Let $\, \x: M \hook \ES^3 \subset \E$
be an  algebraic minimal surface of degree 1.
Then the structure functions satisfy
\beq
h_{\pm 2} = 0. \n
\eeq
The minimal surface  is congruent to a part of the totally geodesic sphere.
\end{theorem}
\noi
Proof of theorem is presented below in 3 steps.

Let us write a real,  homogeneous degree $\,1$ polynomial $\, F \in  S^1(\E)$
in terms of the frame $\, e^{\C}$.
\beq
F  =  p_0\,e_0 + p_1\,E_1+p_{-1}\,E_{-1} + p_3 \, e_3 \in S^1(\E).  \n
\eeq
The equation $\, d F = 0$ is equivalent to the following set of  structure equations
for the polynomial coefficients.
\begin{align}\label{32struct}
&d  p_{{0}}   -p_{{-1}}\omega-p_{{1}}\bomega=0,   \\
&d   p_{{3}}   +p_{{1}}h_{{2}}\omega+p_{{-1}}h_{{-2}}\bomega=0, \n \\
&d   p_{{1}}  - \im \,  p_{{1}} \rho
+\frac{1}{2}\,p_{{0}}\omega
-\frac{1}{2} \,p_{{3}}h_{{-2}}\bomega=0, \n \\
&d   p_{{-1}}   + \im \,  p_{{-1}} \rho
+\frac{1}{2}\,p_{{0}}\bomega
-\frac{1}{2} \,p_{{3}}h_{{2}}\omega=0. \n
\end{align}
\noi
Differentiating the coefficients of $\, F$ would mean applying this structure equation from now on.

\texttt{Step 1}.
Assume the initial condition,
\beq
 p_0 = 0. \n
\eeq
By the metric duality,
this is equivalent to that  $\, F$ vanishes on the minimal surface $\, \x$.

\texttt{Step 2}.
Differentiating $\, p_0=0$, one gets
\beq
 p_{\pm 1}=0. \n
\eeq
At this stage,  the polynomial is reduced to  $\, F = p_3\, e_3$.

\texttt{Step 3}.
Differentiating $\, p_{\pm 1}=0$, one gets
\beq
p_3   h_{\mp 2}=0. \n
\eeq

$\bullet$
If  the structure  functions  $\, h_{\pm 2}$ do  not vanish identically,
$\, p_3 =0$  and the minimal surface does not support
a nonzero degree 1 polynomial vanishing on  the surface.

$\bullet$
If  $\,  h_{\pm 2} = 0$ identically,
the minimal surface is congruent to a part of  the totally geodesic sphere,
Example \ref{21sphere}.
The structure equation \eqref{32struct}
is reduced to   $\, d p_3 =0$, which is compatible.

Fix a point $\, \x_0$ on the minimal surface.
Choose an orthonormal  coordinate $\, \{ \, x_0, \, x_1, \, x_2, \, x_3 \, \}$ of $\, \E$ so that
one has the following identification at $\, \x_0$  by the metric duality.
\begin{align}
e_0  &= x_0, \n \\
E_{\pm 1} &=  \frac{1}{2} ( x_1 \mp  \im \, x_2), \n \\
e_3  &= x_3. \n
\end{align}
Up to scale, one may assume $\, p_{3} = 1$.
The degree 1  polynomial $\, F$ is given by
\beq
F = x_3. \n  \quad \sq
\eeq

\subsection{Degree  2  algebraic minimal surfaces}\label{sec33}
In this section,
we give a  characterization of degree 2  algebraic minimal surfaces.
\begin{theorem}\label{33thm}
Let $\, \x: M \hook \ES^3 \subset \E$ be an  algebraic minimal surface of degree 2.
Then the structure   functions  satisfy
\begin{align}
1 - h_2 h_{-2} &= 0, \n  \\
  h_{\pm 3} &=0. \n
\end{align}
The minimal surface  is congruent to a part of  Clifford torus.
For an orthonormal  coordinate $\, \{ \, x_0, \, x_1, \, x_2, \, x_3 \, \}$
of $\, \E$,
the surface  is defined by the quadratic polynomial
\beq
F = x_0 x_3 - x_1 x_2. \n
\eeq
\end{theorem}
\noi
Proof of theorem is presented below in 5 steps.

Let us write a real,   homogeneous degree $\,2$ polynomial $\, F \in  S^2(\E)$
in terms of the frame $\, e^{\C}$.
\begin{align}
F  &= p_{{0,0}}e_{{0}}^{2}+2\,p_{{0,3}}e_{{0}}e_{{3}} +p_{{3,3}}e_{{3}}^{2}\n \\
    &\quad+2\,r_{{0,1}}e_{{0}}E_{{1}}+2\,r_{{0,-1}}e_{{0}}E_{{-1}}
    +2\,r_{{3,1}}e_{{3}}E_{{1}}+2\,r_{{3,-1}}e_{{3}}E_{{-1}} \n \\
    &\quad+4\,r_{{1,1}}E_{{1}}^{2}+8\,r_{{1,-1}}E_{{1}}E_{{-1}}+4\,r_{{-1,-1}}E_{{-1}}^{2}
    \in S^2(\E). \n
\end{align}
We implicitly assume the appropriate conjugation relations
among the coefficients
so  that $\, F$ is real, i.e.,  $\, r_{0,-1} = \overline{r_{0,1}}$, etc.
The equation $\, d F = 0$ is equivalent to the following set of
structure equations for the polynomial coefficients.
\begin{align}\label{33dF}
&d   p_{{0,0}}   - r_{{0,-1}}\omega -r_{{0,1}}\bomega=0,   \\
&d  p_{{0,3}}
  + \frac{1}{2}\, ( -r_{{3,-1}}+h_{{2}}r_{{0,1}} )   \omega
  + \frac{1}{2} \,  ( -r_{{3,1}}+h_{{-2}}r_{{0,-1}} ) \bomega=0, \n \\
&d  p_{{3,3}}   +h_{{2}}r_{{3,1}}\omega +h_{{-2}}r_{{3,-1}}\bomega=0, \n \\
&d   r_{{0,1}}  - \im \,    r_{{0,1}} \rho
+    (  \,p_{{0,0}}-2\,r_{{1,-1}} ) \omega
+   ( - \,h_{{-2}}p_{{0,3}}-2\,r_{{1,1}} ) \bomega=0, \n \\
 &d   r_{{0,-1}}  + \im  \,  r_{{0,-1}}\rho
+   (  \,p_{{0,0}}-2\,r_{{1,-1}} ) \bomega
+  ( - \,h_{{2}}p_{{0,3}}-2\,r_{{-1,-1}}  ) \omega=0, \n \\
&d   r_{{3,1}}  - \im \, r_{{3,1}}\rho
+    (  \,p_{{0,3}}+2\,h_{{2}}r_{{1,1}}  ) \omega
+ ,  ( - \,h_{{-2}}p_{{3,3}}+2\,h_{{-2}}r_{{1,-1}} ) \bomega=0, \n \\
&d   r_{{3,-1}}   + \im \, r_{{3,-1}}\rho
+    (  \,p_{{0,3}}+2\,h_{{-2}}r_{{-1,-1}} ) \bomega
+   ( - \,h_{{2}}p_{{3,3}}+2\,h_{{2}}r_{{1,-1}}) \omega=0, \n \\
&d  r_{{1,1}}  -2\,\im \,  r_{{1,1}}\rho
+ \frac{1}{2}\,r_{{0,1}}\omega
-\frac{1}{2} \,h_{{-2}}r_{{3,1}}\bomega=0, \n \\
&d r_{{1,-1}}
+\frac{1}{4}\, (  \,r_{{0,-1}}- \,h_{{2}}r_{{3,1}}  ) \omega
+\frac{1}{4}\,  (  \,r_{{0,1}}- \,h_{{-2}}r_{{3,-1}}  ) \bomega=0, \n \\
&d  r_{{-1,-1}} +2\, \im \, r_{{-1,-1}}\rho
+\frac{1}{2}\,r_{{0,-1}}\bomega
 -\frac{1}{2}\,h_{{2}}r_{{3,-1}}\omega=0. \n
\end{align}
\noi
Differentiating the coefficients of $\, F$ would mean applying this structure equation from now on.

\texttt{Step 1}.
Assume the initial condition,
\beq
 Eq_{1,1}: \;  p_{0,0} = 0. \n
\eeq
By the metric duality,
this is equivalent to that
the quadratic polynomial $\, F$  vanishes on the minimal surface $\, \x$.

\texttt{Step 2}.
Differentiating $\, p_{0,0}=0$, one gets
\begin{align}
Eq_{2,1}:  \;&  r_{0,1}   =0,  \n  \\
Eq_{2,2}:  \;&r_{0,-1}  =0. \n
\end{align}

\texttt{Step 3}.
Differentiating $\, r_{0, \pm 1}=0$, one gets
\begin{align}\label{Eq33}
Eq_{3,1}: \;&   h_{{-2}} p_{{0,3}}+2\,r_{{1,1}}  =0,   \\
Eq_{3,2}: \;&  r_{1,-1} =0, \n \\
Eq_{3,3}: \;&  h_{{2}} p_{{0,3}}+2\,r_{{-1,-1}} =0. \n
\end{align}
One may solve for $\, \{ \, r_{1,1}, \, r_{1, -1}, \, r_{-1,-1} \, \}$ from these equations.

\texttt{Step 4}.
Differentiating  \eqref{Eq33}, one gets
\begin{align}\label{Eq34}
Eq_{4,1}: \;&   h_2 \, r_{3,1} =0,   \\
Eq_{4,2}: \;&   h_{-2} \, r_{3,-1} =0. \n
\end{align}
Since the minimal  surface has degree 2,  $\, h_{\pm 2}$ do   not vanish identically,
and this implies $\, r_{3, \, \pm 1} = 0$.

\texttt{Step 5}.
Differentiating  \eqref{Eq34}, one finally gets
\begin{align}\label{Eq35}
Eq_{5,1}: \; &  h_{2} p_{3,3} =0, \n \\
Eq_{5,2}: \; & h_{-2} p_{3,3} =0, \n \\
Eq_{5,3}: \; & (1 - h_2  h_{-2}) \,  p_{0,3} =0. \n
\end{align}

As in \texttt{Step 4},  $\, h_{\pm 2}$ does not vanish identically, and $\, p_{3,3}=0$.
At this stage,  the quadratic  polynomial $\, F$ is reduced to
$\, F =2\,p_{{0,3}}  ( -e_{{0}}e_{{3}}+h_{{-2}} E_{{1}}^{2} +h_{{2}} E_{{-1}}^{2}  )$.

$\bullet$
If  the curvature $\, (1 - h_2  h_{-2})$ does not vanish identically,
$\, p_{{0,3}} =0$  and the minimal surface does not support
a nonzero degree 2 polynomial  vanishing on the surface.

$\bullet$
If  $\,   1 - h_2  h_{-2}  = 0$ identically,
the minimal surface is congruent to Clifford torus,
Example \ref{211torus}.
The structure equation \eqref{33dF} is reduced to
 $\, d p_{{0,3}} =0$, which is compatible.
By the existence  and uniqueness theorem of  ODE,
there  exists up to scale a unique nonzero quadratic  polynomial
that  vanishes on  Clifford  torus.

Fix a point $\, \x_0$ on Clifford  torus.
From the compatibility equation $\,   1 - h_2  h_{-2} = 0$,
one may adapt the frame at $\, \x_0$ so that $\, h_{\pm 2} = 1$.
Choose an orthonormal  coordinate $\, \{ \, x_0, \, x_1, \, x_2, \, x_3 \, \}$ of $\, \E$ so that
one has the following identification at $\, \x_0$ by the metric duality.
\begin{align}
e_0  &= x_0, \n \\
E_{   1} &=  \frac{1}{2} ( x_1 - \im \,  x_2) \, \exp{(\im \frac{\pi}{4})}, \n \\
E_{- 1} &=  \frac{1}{2} ( x_1 + \im \, x_2)  \, \exp{(-\im \frac{\pi}{4})}, \n \\
e_3  &=  x_3. \n
\end{align}
Up to scale,  one may assume $\, p_{0,3} = -\frac{1}{2}$.
The degree 2  polynomial $\, F$ is given by
\beq
F = x_0 x_3 - x_1 x_2. \n  \quad    \sq
\eeq

Comparing the analysis for degree 1, and degree 2  minimal  surfaces,
it  is  evident that
the analysis for  higher degree algebraic  minimal surfaces
would follow the same path,
but  that
the computational  complexity would increase
due  to possibly   the large number of terms
involving  the higher order  structure functions
$\, h_{\pm 2}, \, h_{\pm 3}, ... \, $.\footnotemark
\, We shall show in the next section that
the required differential analysis is manageable for the degree 3 surfaces,
and
one may recover the  Perdomo's result
essentially by  local analysis.

\footnotetext{
This does not necessarily mean  that
the higher degree algebraic minimal surfaces
are characterized by
a set of  higher order  equations.
But the differential analysis itself
does   require  manipulation of
the higher order terms. }

\section{Degree  3  algebraic minimal surfaces}\label{sec4}
In this section,
we give a  local analytic characterization of     degree  3  algebraic minimal surfaces in $\, \ES^3$.
\begin{theorem}\cite{Pe}\label{4thm}
Let $\, \x: M \hook \ES^3 \subset \E$ be an   algebraic minimal surface of degree 3.
Then the structure functions satisfy
\begin{align}
h_2^3 h_{-3}^2 + h_{-2} ^3 h_{3}^2&= 0, \n  \\
h_3 h_{-3} + 4 \,h_2^2 h_{-2}^2 + 4 \, h_2 h_{-2} &= 10 \,(h_2 h_{-2})^{\frac{3}{2}}. \n
\end{align}
For an orthonormal  coordinate $\, \{ \, x_0, \, x_1, \, x_2, \, x_3 \, \}$ of $\,\, \E$,
the surface  is defined by the cubic polynomial
\beq
F = - 2\,x_{{0}}x_{{1}}x_{{2}}+ x_{{3}}( x_{1}^{2} - x_{2}^{2}). \n
\eeq

a)
$\, \x$ is  the conjugate surface  of
a principally bi-planar minimal surface.
It has the Killing nullity 1,
and there  exists up to scale a unique Killing  vector field of  $\, \SO_4$
that is tangent to $\, \x$.

b)
The curvature of the minimal surface
takes values in the closed interval $\, [\,  -3, \, \frac{3}{4} \, ]$.
\end{theorem}
\noi
The  minimal surface defined by the above cubic polynomial
is one of the infinite sequence of  algebraic minimal tori constructed  by Lawson, \cite[p350]{La}.

\begin{remark}
Given the degree 3 algebraic minimal torus,
one may ask if the conjugate surface,
which is principally bi-planar, is also algebraic.
An analysis indicates that
this conjugate minimal surface
does not close up to become a torus.
\end{remark}

Perdomo first gave the characterization of  degree 3 algebraic minimal surfaces in $\, \ES^3$, \cite{Pe}.
One of the main idea of his analysis is   that
such a  minimal surface necessarily contains  a great circle,
and that the gradient of the  defining cubic polynomial
vanishes at  a point on the great circle.
This puts the cubic polynomial into a special  normal form.
By applying the minimal surface equation \eqref{1poly},
the characterization  is reduced essentially to solving  a set of  algebraic equations
among the constant  polynomial coefficients.

The analysis for  the degree 3 case proceeds
similarly as for the degree 1, and    the  degree 2  cases
treated in the previous section.
On the other hand,
a cubic polynomial on $\, \E=\R^4$ has 20 coefficients.
The computational complexity increases as degree increases,
and
one  has to take higher order  derivatives
in order to access the compatibility equations for  a minimal surface to be algebraic.
For the degree 3 case,
the analysis  requires  differentiating  six times.

Proof of  Theorem \ref{4thm}  consists of   two   parts.
In the first part,
a preliminary analysis is  done
to reduce the number of independent polynomial coefficients from 20 to 6,  Section \ref{sec41}.
After differentiating four times,
the analysis divides into two  cases
depending on whether a certain structure invariant of the minimal surface vanishes or not.
In the second part,
we carry out the differential analysis for each case
under the appropriate assumptions on the structure functions,
Section \ref{sec42},  Section \ref{sec44}.

For the rest of the paper,
we assume $\, h_{\pm 2} \ne  0, \, 1 - h_2 h_{- 2} \ne 0$,
and  $\, h_{\pm 3} \ne  0$
to exclude   degree 1, and   degree 2 algebraic  minimal surfaces.

\subsection{Preliminary analysis}\label{sec41}

Let us write a real,  homogeneous degree $\,3$ polynomial $\, F \in  S^3(\E)$
in terms of the frame $\, e^{\C}$, \eqref{21ecframe}.
\begin{align}\label{41F}
F  =
&\;\;
p_{0}e_{0}^{3}+p_{1}e_{0}^{2}e_{3}+p_{2}e_{0}e_{3}^{2}+p_{{3}}e_{{3}}^{3}   \\
&
+8\,q_{0}E_{1}^{3}+8\,q_{1}E_{1}^{2}E_{-1}+8\,q_{2}E_{1}E_{-1}^{2}+8\,q_{{3}}E_{-1}^{3}\n \\
&
+2\,r_{{0,1}}e_{{0}}^{2}E_{{1}}+2\,r_{{0,-1}}e_{{0}}^{2}E_{{-1}}\n \\
&
+2\,r_{{1,1}}e_{{0}}e_{{3}}E_{{1}}+2\,r_{{1,-1}}e_{{0}}e_{{3}}E_{{-1}}\n \\
&
+2\,r_{{2,1}}e_{{3}}^{2}E_{{1}}+2\,r_{{2,-1}}e_{{3}}^{2}E_{{-1}}\n \\
&
+4\,r_{0,2}e_{0}E_{1}^{2}+4\,r_{0,0}e_{0}E_{1}E_{-1}+4\,r_{0,-2}e_{0}E_{-1}^{2}\n \\
&
+4\,r_{{1,2}}e_{{3}}E_{{1}}^{2}+4\,r_{{1,0}}e_{{3}}E_{{1}}E_{{-1}}
+4\,r_{{1,-2}}e_{{3}}E_{{-1}}^{2}. \n
\end{align}
We implicitly assume the appropriate conjugation relations
among the coefficients
so  that $\, F$ is real, i.e.,  $\, r_{0,-1} = \overline{r_{0,1}}$, etc.

$\, F \in  S^3(\E)$ is a constant cubic polynomial.
Differentiating    \eqref{41F}    by Leibniz  rule,
the equation $\, d F = 0$ is equivalent to the following set of
structure equations for the polynomial coefficients.

\begin{align}\label{41struct}
&
d   p_{{0}} -r_{{0,-1}}\omega -r_{{0,1}}\bomega=0,   \\
&
d  p_{{1}}   +\left( -r_{{1,-1}}+h_{{2}}r_{{0,1}} \right) \omega
  + \left( -r_{{1,1}}+ h_{{-2}}r_{{0,-1}} \right) \bomega=0, \n \\
&
d  p_{{2}} +\left( -r_{{2,-1}}+h_{{2}}r_{{1,1}} \right) \omega
+ \left( -r_{{2,1}}+h_{{-2}}r_{{1,-1}} \right) \bomega=0, \n \\
&
d  p_{{3}} +h_{{2}}r_{{2,1}}\omega +h_{{-2}}r_{{2,-1}}\bomega=0, \n  \\
&
d   q_{{0}}
-3\,\im \, q_{{0}}\rho+\frac{1}{2}  r_{{0,2}}\omega
-\frac{1}{2} h_{{-2}}r_{{1,2}}\bomega=0, \n  \\
&
 d   q_{{1}}
-\im \, q_{{1}}\rho+\frac{1}{2} \left(  \,r_{{0,0}}- \,h_{{2}}r_{{1,2}} \right)\omega
+\frac{1}{2}  \left(  \,r_{{0,2}}- \,h_{{-2}}r_{{1,0}} \right) \bomega=0, \n  \\
&
 d   q_{{2}}
+\im \, q_{{2}}\rho
+\frac{1}{2}  \left(  \,r_{{0,-2}}- \,h_{{2}}r_{{1,0}} \right) \omega
+\frac{1}{2}  \left(  \,r_{{0,0}}- \,h_{{-2}}r_{{1,-2}} \right) \bomega=0, \n  \\
&
 d  q_{{3}}
+3\,\im \, q_{{3}}\rho
-\frac{1}{2} h_{{2}}r_{{1,-2}}\omega+\frac{1}{2} r_{{0,-2}}\bomega=0, \n
\end{align}
\begin{align}
&
d   r_{{0,1}}  -  \im \,   r_{{0,1}}\rho
+\frac{1}{2}  \left( -2\,r_{{0,0}}+3\,p_{{0}} \right) \omega
+\frac{1}{2}  \left( -4\,r_{{0,2}}  - h_{{-2}}p_{{1}}\right) \bomega=0, \n  \\
&
 d   r_{{0,-1}}  + \im \,  r_{{0,-1}}\rho
+\frac{1}{2}  \left( -4\,r_{{0,-2}}-h_{{2}}p_{{1}} \right)\omega
+\frac{1}{2}  \left( -2\,r_{{0,0}}+3\,p_{{0}} \right) \bomega=0, \n  \\
&
 d   r_{{1,1}}   - \im \,  r_{{1,1}}\rho
+ \, \left( - \,r_{{1,0}}+ \,p_{{1}}+2\,h_{{2}}r_{{0,2}} \right) \omega
+ \, \left( - 2 \,r_{{1,2}} - \,h_{{-2}}p_{{2}}+ \,h_{{-2}}r_{{0,0}} \right) \bomega=0, \n  \\
&
d   r_{{1,-1}}  + \im \,  r_{{1,-1}}\rho
+ \, \left(  \,-2\,r_{{1,-2}}- \,h_{{2}}p_{{2}} +h_{{2}}r_{{0,0}}\right) \omega
+ \, \left( - \,r_{{1,0}}+ \,p_{{1}} +2\,h_{{-2}}r_{{0,-2}} \right) \bomega=0, \n  \\
&
 d  r_{{2,1}} - \im \,  r_{{2,1}}\rho
+\frac{1}{2}\, \left( p_{{2}}+4\,h_{{2}}r_{{1,2}} \right) \omega
+\frac{1}{2}\, \left( 2\,h_{{-2}}r_{{1,0}}-3\,h_{{-2}}p_{{3}} \right) \bomega=0, \n  \\
&
 d   r_{{2,-1}}  + \im \,  r_{{2,-1}}\rho
+\frac{1}{2}\, \left( 2\,h_{{2}}r_{{1,0}}-3\,h_{{2}}p_{{3}} \right) \omega
+\frac{1}{2}\, \left( p_{{2}}+ 4\,h_{{-2}}r_{{1,-2}} \right) \bomega=0, \n
\end{align}
\begin{align}
&
 d  r_{{0,2}}   -2\,\im \,   r_{{0,2}}\rho
+ \, \left( - \,q_{{1}}+ \,r_{{0,1}} \right) \omega
+\frac{1}{2}\, \left( -6 \,q_{{0}}- \,h_{{-2}}r_{{1,1}} \right) \bomega=0, \n  \\
&
 d  r_{{0,0}}
+\frac{1}{2}\, \left( -4\,q_{{2}}- \,h_{{2}}r_{{1,1}}+2\,r_{{0,-1}} \right) \omega
+\frac{1}{2}\, \left(- 4 \,q_{{1}} - \,h_{{-2}}r_{{1,-1}}+2\,r_{{0,1}} \right)\bomega=0, \n  \\
&
d   r_{{0,-2}} +2\,\im \,  r_{{0,-2}}\rho
 +\frac{1}{2}\, \left( -6\,q_{{3}}  - \,h_{{2}}r_{{1,-1}} \right) \omega
 + \, \left( - \,q_{{2}}+ \,r_{{0,-1}}\right) \bomega=0, \n  \\
&
d  r_{{1,2}}   -2\,\im \,  r_{{1,2}}\rho
+\frac{1}{2}\, \left(  \,r_{{1,1}}+6\,h_{{2}}q_{{0}} \right) \omega
+ \, \left( - \,h_{{-2}}r_{{2,1}}+ \,h_{{-2}}q_{{1}} \right) \bomega=0, \n  \\
&
 d   r_{{1,0}}
+\frac{1}{2}\, \left(  \,r_{{1,-1}}+4\,h_{{2}}q_{{1}}-2\,h_{{2}}r_{{2,1}} \right) \omega
+\frac{1}{2}\, \left(  \,r_{{1,1}}+4\,h_{{-2}}q_{{2}} -2\,h_{{-2}}r_{{2,-1}}\right) \bomega=0, \n  \\
&
 d   r_{{1,-2}}  +2\,\im \,   r_{{1,-2}}\rho
+ \, \left( - \,h_{{2}}r_{{2,-1}}+ \,h_{{2}}q_{{2}} \right) \omega
+\frac{1}{2}\, \left(  \,r_{{1,-1}}+6\,h_{{-2}}q_{{3}} \right) \bomega=0. \n
\end{align}
\noi
Differentiating the coefficients of $\, F$ would mean applying this structure equation from now on.

The preliminary analysis consists of  the following 5 steps.
The condition that $\, F$ vanishes on the minimal surface
serves as  the initial condition  for the over-determined PDE analysis.
By successively applying the above structure equation,
we reduce the number of independent polynomial coefficients to $\, 20 - 14 = 6$.
The reduction process stops  at \texttt{Step 5},
where  the analysis divides into two  cases.

\texttt{Step 1}.
Assume the initial condition
\beq
 Eq_{1,1}: \;  p_{0} = 0. \n
\eeq
By the metric duality,
this equivalent to that
the cubic  polynomial $\, F$  vanishes on the minimal surface $\, \x$.

\texttt{Step 2}.
Differentiating $\, p_{0 }=0$, one gets
\begin{align}
Eq_{2,1}:  \;&  r_{0,1}   =0,  \n  \\
Eq_{2,2}:  \;&r_{0,-1}  =0. \n
\end{align}

\texttt{Step 3}.
Differentiating $\, r_{0, \pm 1}=0$, one gets
\begin{align}\label{Eq43}
Eq_{3,1}: \;& \frac{1}{2} \,h_{{-2}}p_{{1}}+2\,r_{{0,2}}  =0,   \\
Eq_{3,2}: \;&  r_{0,0} =0, \n \\
Eq_{3,3}: \;& \frac{1}{2}\,h_{{2}}p_{{1}}+2\,r_{{0,-2}} =0. \n
\end{align}
One may solve for $\, \{ \, r_{0,2}, r_{0,0}, \, \, r_{0,-2} \, \}$ from these equations.

\texttt{Step 4}.
Differentiating  \eqref{Eq43}, one gets
\begin{align}\label{Eq44}
Eq_{4,1}: \;&   2\,q_{{2}}+ \frac{1}{2}\,h_{{2}}r_{{1,1}} =0,    \\
Eq_{4,2}: \;&   2\,q_{{1}}+ \frac{1}{2}\,h_{{-2}}r_{{1,-1}} =0,  \n   \\
Eq_{4,3}: \;&  6\,q_{{0}} + \frac{1}{2}\,p_{{1}}h_{{-3}}+\frac{3}{2}\,h_{{-2}}r_{{1,1}}  =0,  \n   \\
Eq_{4,4}: \;&  6\,q_{{3}} +  \frac{1}{2}\,p_{{1}}h_{{3}}+\frac{3}{2}\,h_{{2}}r_{{1,-1}}  =0. \n
\end{align}
One may solve for $\, \{ \, q_0, \, q_1,\, q_2, \, q_3 \, \}$ from these equations.

\texttt{Step 5}.
Differentiating  \eqref{Eq44}, one   gets
\begin{align}\label{Eq45}
Eq_{5,1}: \; &
2 \,h_{{3}} r_{{1,1}}+6 \,h_{{2}}r_{{1,0}}
+(- \,h_{{2}} +  \,h_{2}^{2}h_{{-2}}) p_{{1}} =0,   \\
Eq_{5,2}: \; &
h_{{-2}}r_{{1,-2}}+h_{{2}}r_{{1,2}}+\frac{1}{2}\,h_{{2}}h_{{-2}}p_{{2}}=0, \n \\
Eq_{5,3}: \; &
2 \,h_{{-3}} r_{{1,-1}}+6 \,h_{{-2}}r_{{1,0}}
+(- \,h_{{-2}} +  \,h_{-2}^{2}h_{{2}}) p_{{1}} =0, \n \\
Eq_{5,4}: \; &
12\,h_{{-3}}h_{2}^{2}r_{{1,0}}
+6\,h_{{3}}h_{-2}^{2}r_{{1,-2}}-6\,h_{{-2}}h_{{3}}h_{{2}}r_{{1,2}}
+(-h_{{3}}h_{{2}}h_{{-4}}   -2\,h_{{-3}}h_{2}^{2}
+2\,h_{{-3}}h_{2}^{3}h_{{-2}} )  p_{{1}}=0, \n \\
Eq_{5,5}: \; &
12 \left(  \,h_{3}^{2}h_{-2}^{3}+ \,h_{-3}^{2}h_{2}^{3} \right) r_{{1,0}}
+
 ( -h_{{3}}h_{-2}^{2}h_{{-3}}h_{{4}}-2\,h_{3}^{2}h_{-2}^{3}
+2\,h_{3}^{2}h_{-2}^{4}h_{{2}}-h_{{-3}}h_{2}^{2}h_{{3}}h_{{-4}}\n \\
&
-2\,h_{-3}^{2}h_{2}^{3}+2\,h_{-3}^{2}h_{2}^{4}h_{{-2}} ) p_{{1}}
=0. \n
\end{align}
Since  we are assuming that
the algebraic  minimal  surface has degree 3,
both $\, h_{\pm 2}, \,  h_{\pm 3}$ do  not vanish identically,
otherwise the algebraic  minimal  surface would have  degree 1, or 2
by  Theorem \ref{32thm}, Theorem \ref{33thm}.
One may thus  solve for $\, \{ \, r_{1,1}, \, p_2, \, r_{1,-1}, r_{1,2} \, \}$ from
$\, \{ \, Eq_{5,1},  \,Eq_{5,2},  \,Eq_{5,3},  \,Eq_{5,4} \, \}$.
Note that $\, Eq_{5,5}$ is equivalent to $\, \overline{Eq_{5,4}}$
modulo $\, Eq_{5,4}$.
$\sq$

At  this  step,
the analysis  is  divided into  the following  two  cases.
Set
\begin{align}\label{41Delta}
\Delta_3^{+} &= h_{{-2}}^{3} h_{{3}}^{2} +  \,h_{{2}}^{3}h_{{-3}}^{2},  \\
\Delta_4 &= h_2 h_{-3}^2 h_{4} -  h_{-2}  h_{3}^2 h_{-4}. \n
\end{align}

$\bullet$ Case $\, \Delta_3^+=0$, Section \ref{sec42}:
It will be shown that
the degree 3 algebraic minimal surface in Theorem \ref{4thm}
belongs to  this case.
The minimal surfaces with $\, \Delta_3^+=0$ are
the conjugate surfaces of the \emph{principally bi-planar} minimal surfaces.
A minimal surface in $\, \ES^3$ is principally bi-planar
when each of its principal curves is planar, \cite{Yam}.
The principally bi-planar minimal surfaces are  characterized  by  the equation
$\, \Delta_3^{-} = h_{{-2}}^{3} h_{{3}}^{2} -  \,h_{{2}}^{3}h_{{-3}}^{2}=0$.\footnotemark
\footnotetext{
One may verify this by a straightforward moving frame computation.
We  omit the details.}
\,
Both of the $\,  \Delta_3^{\pm}$-null  minimal surfaces
belong to the wider class of
the cohomogeneity 1 minimal surfaces, \cite[p32]{HsL}.\footnotemark

\footnotetext{
One may verify by direct computation that
the  cohomogeneity 1 minimal surfaces in $\, \ES^3$
are characterized by the pair of  fourth order equations
$\, \Delta^+_4  =
2\, h_2 h^2_{-3} h_4  - h_3 h_{-3} ( 3\,h_3 h_{-3} + 2\,h_2 h_{-2} \, K) =0, \,
\Delta^-_{4}  = \overline{ \Delta^+_4}=0$.}

$\bullet$ Case $\,\Delta_3^{+} \not \equiv  0, \,  \Delta_4 = 0\,$, Section \ref{sec44}:
The analysis shows that
for a minimal surface with $\, \Delta_3^{+}  \not \equiv 0$,
the fourth order equation  $\, \Delta_4 = 0$
is  a necessary condition  to be algebraic of degree 3.
But, a further analysis shows that
the resulting structure equation is not compatible,
and
there does not exist any degree 3 algebraic minimal surfaces in this  case.

In the following two sections,
we present the differential analysis  for each case.

\subsection{Case $\,    h_{{-2}}^{3} h_{{3}}^{2} +  \,h_{{2}}^{3}h_{{-3}}^{2}=0$}
\label{sec42}
In this section,
we show that, up to motion by $\, \SO_4$,
there exists a unique degree 3 algebraic minimal surface in $\, \ES^3$
that satisfies
$\, \Delta_3^{+}  =  h_{{-2}}^{3} h_{{3}}^{2} +  \,h_{{2}}^{3}h_{{-3}}^{2}=0$.

We first determine the structure equation for the minimal surfaces
with $\, \Delta_3^{+}=0$,  Section \ref{sec421}.
A first integral $\, \JAI$  is defined, \eqref{421JAI},
and  it follows that
there exists locally a one parameter family of   $\, \Delta_3^{+}$-null surfaces.
We then continue the analysis of Section \ref{sec41}
and show that
exactly one of the  $\, \Delta_3^{+}$-null surfaces is algebraic of degree 3, Section \ref{sec422}.

\subsubsection{Structure equation}\label{sec421}
Let $\, \x: M \hook \ES^3$ be a  $\, \Delta_3^+$-null surface.
Differentiating $\, \Delta_3^{+}=0$,
one may solve for $\, h_{\pm 4}$ and get
\begin{align}\label{421h4}
h_{{4}}&=
\frac{1}{2} \,{\frac {h_{2}^{2}h_{{-3}}
( -2\,h_{{2}}h_{{-2}}+2\,h_{2}^{2}h_{-2}^{2}-3\,h_{{3}}h_{{-3}}) }
{h_{{3}}h_{-2}^{3}}} ,    \\
h_{{-4}}&=\overline{h_4}. \n
\end{align}
Here we assume $\, h_{\pm 2}, \, h_{\pm 3}$ are nonzero.
A  direct computation shows that
the   structure  equation is compatible with this relation, i.e.,
$\,  d^2h_{\pm 3}=0$ is an identity.

Set
\beq\label{421JAI}
\JAI =\frac{1}{8} {\frac { \left( h_{{3}}h_{{-3}}+4\,h_{{2}}h_{{-2}}
+4\,h_{2}^{2}h_{-2}^{2} \right) }{(h_2h_{-2})^{\frac{3}{2}}}}, \; \; \JAI >0.
\eeq
Then
$\, d \JAI=0$,
and  $\,  \JAI$ is a  first integral  for the  $\, \Delta_3^+$-null surfaces.

We claim that $\, \JAI  \in  ( 1, \, \infty)$, and that $\, h_2$ is nowhere zero.
Let us denote  $\, | h_2 | = a, \, |h_3|=b $.
The equation \eqref{421JAI}  is written as
\beq\label{baJ}
b^2+4\,a^2( 1 +a^2-2\, \JAI \, a) = 0.
\eeq
This implies that
$\, \JAI \geq 1$, and that
either $\, a  \equiv  0$, which is excluded, or
$\, a$ takes values in the closed interval $\, [ \JAI - \sqrt{\JAI^2 - 1}, \, \JAI + \sqrt{\JAI^2 - 1}]$.
When $\, \JAI =1$, one has that $\, a  \equiv   1, \, b  \equiv   0$, which is also excluded.

The structure equations
\eqref{421h4} and \eqref{421JAI} will be used implicitly
for the  analysis in  the next  subsection.


\subsubsection{Differential analysis}\label{sec422}
We now continue the analysis of Section \ref{sec41}.
Due to their lengths,
the exact  expressions  for
$\, Eq_{6,1}, \, Eq_{6,2}, Eq_{6,3}, \, Eq_{6,4}; \, Eq_7$ below
will be postponed to Appendix.

\texttt{Step  5'}.
Assume $\, \Delta_3^+=0$.
Then  \eqref{421h4}  implies   that  $\, Eq_{5,5}=0$.

\texttt{Step 6}.
Differentiating   \eqref{Eq45} and equating modulo $\, \Delta_3^+=0$,
one gets a set of  four  independent equations
$\, \{ \, Eq_{6,1}, \, Eq_{6,2}, Eq_{6,3}, \, Eq_{6,4} \, \}$(see Appendix).
One may solve for
$\, \{ \,  r_{1,0}, \, r_{2,1}, \, r_{2,-1} r_{1,-2} \, \}$
from these equations.

\texttt{Step 7}.
Differentiating  $\, Eq_{6,1}$,
one gets $\, Eq_7$(see Appendix).
One may solve for
$\, p_3$ from this equation.

At this step,  $\, p_1$ is the only remaining independent  polynomial coefficient,
and it satisfies the structure equation,  \eqref{41struct},  of the form
\beq
dp_1 \equiv 0, \mod \; \; p_1. \n
\eeq
From the uniqueness theorem of ODE,
if $\, p_1$ vanishes at a point of the minimal surface,
it vanishes identically, which implies
the cubic polynomial $\, F$ vanishes.
We therefore assume $\, p_1$ is nowhere zero from now on.

\texttt{Step 8}.
Differentiating the remaining equations $\, \{ \, Eq_{6,2}, Eq_{6,3}, \, Eq_{6,4}; \, Eq_7 \, \}$,
one gets a single compatibility equation,  up to scale by non-identically zero terms;
\beq
(8 \JAI+10)(8 \JAI -10)=0. \n
\eeq
Here $\, \JAI$ is the first integral   \eqref{421JAI}.
Since $\, \JAI > 1$, one must have
\beq\label{J54}
\JAI = \frac{5}{4}.
\eeq
Moreover, a direct computation shows that $\, d^2 p_1=0$ is an identity
with this relation.

By the existence  and uniqueness theorem of  ODE,
there  exists up to scale a unique nonzero cubic polynomial
that  vanishes on the $\, \Delta_3^+$-null minimal surface
with the first integral $\, \JAI = \frac{5}{4}$.

Fix a point $\, \x_0$  on the minimal surface.
Choose an orthonormal  coordinate $\, \{ \, x_0, \, x_1, \, x_2, \, x_3 \, \}$ of $\, \E$ so that
one has the following identification at $\, \x_0$ by the  metric duality.
\begin{align}
e_0  & =   x_1, \n \\
E_{   1} & =     \frac{1}{2} ( x_0 - \im \, x_2) \, \exp{(\im \frac{\pi}{4})}, \n \\
E_{- 1} & =     \frac{1}{2} ( x_0 + \im \, x_2)  \, \exp{(-\im \frac{\pi}{4})}, \n \\
e_3  & =    x_3. \n
\end{align}
Up to scale,  one may assume $\, p_{1} =  \frac{1}{3}$ at $\, \x_0$.
From the analysis of Section \ref{sec421},
$\, |h_2| = a$ takes values in the closed interval $\, [ \, \frac{1}{2}, \,  2 \, ]$.
Evaluating  $\, \lim_{a \to 2} F$ with $\, h_2 = a$,
the degree  3   polynomial $\, F$ is given by
\beq\label{Fformula}
F = - 2\,x_{{0}}x_{{1}}x_{{2}}+ x_{{3}}( x_{1}^{2} - x_{2}^{2}).
\eeq

\subsubsection{Proof of Theorem \ref{4thm}}\label{sec423}

\

a)
From \eqref{Fformula},
set  $\, z_1 = x_3 + \im \, x_0, \, z_2= x_1 + \im \, x_2$.
Then $\, F = Re( z_1  z_2^2)$.
It is clear that $\, F$ is invariant under a subgroup $\, \SO_2 \subset \SO_4$.
It   is known that
a minimal surface in $\, \ES^3$ with Killing nullity $\, \geq 2$
is  either the totally geodesic sphere or Clifford torus, \cite{HsL}.

b)
By definition of the first integral $\, \JAI$ in  \eqref{baJ},
$\, | \, h_2   | = a$ takes values in the closed interval $\, [ \, \frac{1}{2}, \, 2 \, ]$.
The curvature is  given by  $\, K  = 1 - a^2$.
$\sq$

\subsection{Case
 $\,  h_{{-2}}^{3} h_{{3}}^{2} +  \,h_{{2}}^{3}h_{{-3}}^{2} \ne     0$}\label{sec44}

In this section,
we show that
there does not exist a  degree 3 algebraic minimal surface in  $\, \ES^3$
with the property that
$\, \Delta_3^+=h_{{-2}}^{3} h_{{3}}^{2} +  \,h_{{2}}^{3}h_{{-3}}^{2}$
is not identically zero.

The  analysis  will show that under the condition   $\,  \Delta_3^+ \not \equiv 0$,
a  degree 3 algebraic minimal surface
necessarily satisfies a fourth order equation
$\, \Delta_4  = h_2 h_{-3}^2 h_{4} -  h_{-2}  h_{3}^2 h_{-4}=0$.
We  first determine the structure equation for the minimal surfaces
with $\, \Delta_4 =0$,  Section \ref{sec441}.
We then continue the analysis of Section \ref{sec41} and show that
none of the  $\, \Delta_4$-null surfaces is algebraic of degree 3,  Section \ref{sec442}.

\subsubsection{Structure equation}\label{sec441}
Let $\, \x: M \hook \ES^3$ be a  $\, \Delta_4$-null surface.
Differentiating $\, \Delta_4=0$,
one may solve for $\, h_{\pm 5}$ and get
\begin{align}\label{441h5}
h_{{5}}&=
\frac{5\,h_{{-2}}h_{{3}}^{2}h_{{-3}}-7\,h_{{-2}}^{2}h_{3}^{2}h_{{-3}}h_{{2}}
+4\,h_{{4}}h_{{-4}}h_{{3}}h_{{-2}}-4\,h_{{4}}h_{{-3}}h_{{2}}h_{{-2}}
+4\,h_{{4}}h_{{-3}}h_{{2}}^{2}h_{{-2}}^{2}-2\,h_{{4}}h_{-3}^{2}h_{{3}}}
{2 \,h_{-3}^2 h_2},    \\
h_{{-5}}&=\overline{h_5}. \n
\end{align}
Here we assume $\, h_{\pm 2}, \, h_{\pm 3}$ are nonzero.
A  direct computation shows that
the   structure  equation is compatible with this relation, i.e.,
$\,  d^2h_{\pm 4}=0$ is an identity.

The structure equation  \eqref{441h5}  will be used implicitly
for the  analysis in  the next  subsection.

\subsubsection{Differential analysis}\label{sec442}
We now continue the analysis of Section \ref{sec41}.
Due to their lengths,
the exact  expressions  for
$\, Eq_{6,1}, \, Eq_{6,2}, Eq_{6,3}; \, Eq_7$ below
will be postponed to Appendix.

\texttt{Step  5'}.
Assume $\, \Delta_3^+\ne 0$.
One may solve for  $\, r_{1,0}$ from  $\, Eq_{5,5}$.

\texttt{Step 6}.
Differentiating   $\, Eq_{5,1}, \, Eq_{5,2}$ from \eqref{Eq45},
one gets a set of  three  independent equations
$\, \{ \, Eq_{6,1}, \, Eq_{6,2}, Eq_{6,3} \, \}$(see Appendix).
One may solve for
$\, \{ \,  r_{2,1}, \, r_{1,-2}, \, r_{2,-1} \, \}$
from these equations.

\texttt{Step 7}.
Differentiating  $\, Eq_{6,1}$,
one gets a single equation  $\, Eq_7$(see Appendix).
One may solve for $\, p_3$ from this equation.

At this step,  $\, p_1$ is the only remaining independent  polynomial coefficient,
and it satisfies the structure equation, \eqref{41struct},  of the form
\beq
dp_1 \equiv 0, \mod \; \; p_1. \n
\eeq
As in Section \ref{sec422}, we assume $\, p_1$ is nowhere zero from now on.

\texttt{Step 8}.
Differentiating
$\, Eq_{5,4}, \, Eq_{5,5}$ from \eqref{Eq45},
one gets a set of two equations,
which allow one to solve for $\, h_{\pm 5}$(see Appendix).
Note that we have not assumed $\, \Delta_4 = 0$ yet.

At this step, we observe that the coefficient $\, p_3$ is real.
Evaluating $\, p_3  - \overline{p_3} =0$ with the relations obtained so far,
one gets the following compatibility equation;
\begin{align}\label{442Key}
&
 \Delta_4 \,  \Delta_4'   =0, \n
\end{align}
where
\begin{align}
\Delta_4 '=&
 \,  ( -4\,h_{2}^{4}h_{-3}^{2}h_{{-2}}
 +3\,h_{2}^{4}h_{-2}^{2}h_{{-4}}+10\,h_{2}^{3}h_{-3}^{2}
-3\,h_{2}^{3}h_{{-2}}h_{{-4}}+3\,h_{2}^{2}h_{{4}}h_{-2}^{4}
-3\,h_{{2}} h_{{-2}}^{3}h_{{4}}\n \\
&  \;
-4\,h_{{2}}h_{-2}^{4}h_{3}^{2}+10\,h_{-2}^{3}h_{3}^{2}  ) . \n
\end{align}
A short analysis shows that
$\, \Delta_4' $ vanishes only when $\, h_3 \equiv 0$, which is excluded.
Hence a degree 3 algebraic minimal surface with
$\, \Delta_3^+ \not \equiv 0$ must satisfy $\, \Delta_4 = 0$.

We assume the structure equations \eqref{441h5}  from now on.

\texttt{Step 9}.
The cubic polynomial  $\, F$, \eqref{41F},  is  real.
Evaluating $\,F - \overline{F} =0$,
one gets the single   compatibility equation
\beq
( -h_{{-2}}^{3}h_{{3}}^{2}+h_{{2}}^{3}h_{{-3}}^{2})
\Delta_4'' =0, \n
\eeq
where
\beq
\Delta_4'' =3\,h_{{2}}h_{{-3}}h_{4}^{2}
+( 3\,h_{2}^{2}h_{{3}} h_{{-2}}^{2}
-3\,h_{{2}}h_{{3}}h_{{-2}}-4\,h_{{-3}}h_{3}^{2}  ) h_{{4}}
-4\,h_{-2}^{2}h_{{2}}h_{3}^{3}+10\,h_{{-2}}h_{3}^{3}. \n
\eeq

A further differential analysis shows that
$\,  ( -h_{-2}^{3}h_{3}^{2}+h_{2}^{3}h_{-3}^{2})$ vanishes
only when $\, h_3 \equiv 0$, which is excluded
(the differential analysis for this case is a little  evolved,
but straightforward. We shall omit the details).
Hence a degree 3 algebraic minimal surface with
$\, \Delta_3^+ \not \equiv 0$ must also satisfy $\, \Delta_4''  = 0$.

\texttt{Step 10}.
Differentiating   $\,  Eq_{6,2}$ modulo $\, \Delta_4''=0 $,
one gets another  single compatibility equation,   up to scale by non-identically zero terms;
\beq
Eq_{10,1}: \left( 3\,h_{{2}}h_{{-2}}-2 \right) h_{{4}}-4\,h_{3}^{2}h_{{-2}}=0.  \n
\eeq
Comparing $\, Eq_{10,1}$ with $\, \Delta_4'' =0 $,
one gets
\beq
Eq_{10,2}: 39\,h_{2}^{2}h_{-2}^{2}-56\,h_{{2}}h_{{-2}}+20+16\,h_{{3}}h_{{-3}}=0. \n
\eeq
Differentiating this equation again, one gets
\beq
Eq_{10,3}: 93\,h_{2}^{2}h_{-2}^{2}-122\,h_{{2}}h_{{-2}}+40+32\,h_{{3}}h_{{-3}}=0. \n
\eeq
$\, Eq_{10,2}$ and $\, Eq_{10, 3}$ are compatible only when $\, h_{\pm 2}=0$.
$\sq$

\np
\renewcommand{\theequation}{A-\arabic{equation}}
\setcounter{equation}{0}
\section*{Appendix. }
We  record the exact formulae of
 the  long expressions omitted in the main text.
\vsp{1pc}

\noi
A-1. \textbf{Section \ref{sec42}}
\begin{align}
Eq_{6,1}&:
 ( -6\,h_{2}^{3}+6\,h_{2}^{4}h_{{-2}}  ) r_{{1,0}}
-6\,h_{2}^{2}h_{{-2}}h_{{3}} r_{{2,-1}}
+4\,h_{3}^{2}h_{{-2}}r_{{1,-2}} \n \\
&\quad
+( h_{2}^{5}h_{-2}^{2}-h_{2}^{3}h_{{-3}}h_{{3}}
+h_{2}^{3}-2\,h_{2}^{4}h_{{-2}} ) p_{{1}}=0
, \n \\
Eq_{6,2}&:
  ( 36\,h_{2}^{6}h_{-2}^{4}h_{{-3}}h_{{3}}-36\,h_{2}^{
5}h_{-2}^{3}h_{{-3}}h_{{3}}  ) r_{{2,1}}
+ ( 54\,h_{-3}^{3}h_{2}^{5}h_{{-2}}h_{{3}}
+72\,h_{3}^{2}h_{-2}^{5}h_{2}^{3}
\n \\
&\quad
-72\,h_{2}^{4}h_{-2}^{6}h_{3}^{2}
+36\,h_{-3}^{2}h_{2}^{6}h_{-2}^{2}
+48\,h_{{-3}}h_{2}^{2}h_{3}^{3}h_{-2}^{4}
-36\,h_{-3}^{2}h_{2}^{7}h_{-2}^{3} ) r_{{2,-1}}
\n \\
&\quad
+  ( -36\,h_{-3}^{3}h_{2}^{3}h_{{-2}}h_{3}^{2}
-32\,h_{3}^{4}h_{-2}^{4}h_{{-3}}+24\,h_{2}^{5}h_{-3}^{2}h_{-2}^{3}h_{{3}}
\n \\
&\quad
-48\,h_{{2}}h_{3}^{3}h_{-2}^{5}
+48\,h_{3}^{3}h_{-2}^{6}h_{2}^{2}
-24\,h_{2}^{4}h_{-3}^{2}h_{{3}}h_{-2}^{2}  ) r_{{1,-2}}
\n \\
&\quad
+  ( -6\,h_{2}^{8}h_{-3}^{3}h_{-2}^{2}
+8\,h_{3}^{3}h_{-3}^{2}h_{2}^{3}h_{-2}^{3}
+6\,h_{-3}^{3}h_{2}^{7}h_{{-2}}+10\,h_{2}^{3}h_{-2}^{3}h_{{-3}}h_{3}^{2}
 \n \\
&\quad
+9\,h_{2}^{6}h_{-3}^{4}h_{{3}}-2\,h_{-2}^{5}h_{2}^{5}h_{{-3}}h_{3}^{2}
-8\,h_{2}^{4}h_{-2}^{4}h_{{-3}}h_{3}^{2} ) p_{{1}}=0
, \n \\
Eq_{6,3}&:
 ( -36\,h_{2}^{4}h_{-2}^{6}h_{3}^{2}
 +102\,h_{-3}^{3}h_{2}^{5}h_{{-2}}h_{{3}}
 -36\,h_{-3}^{2}h_{2}^{7}h_{-2}^{3}
 \n \\
&\quad
 +36\,h_{3}^{2}h_{-2}^{5}h_{2}^{3}
 +36\,h_{-3}^{2}h_{2}^{6}h_{-2}^{2}
 +48\,h_{{-3}}h_{2}^{2}h_{3}^{3}h_{-2}^{4} t) r_{{2,-1}}
 \n \\
&\quad
+ ( 24\,h_{3}^{3}h_{-2}^{6}h_{2}^{2}
-24\,h_{{2}}h_{3}^{3}h_{-2}^{5}-48\,h_{2}^{4}h_{-3}^{2}h_{{3}}h_{-2}^{2}
\n \\
&\quad
+48\,h_{2}^{5}h_{-3}^{2}h_{-2}^{3}h_{{3}}
-68\,h_{-3}^{3}h_{2}^{3}h_{{-2}}h_{3}^{2}
-32\,h_{3}^{4}h_{-2}^{4}h_{{-3}} ) r_{{1,-2}}
\n \\
&\quad
+ (14\,h_{3}^{3}h_{-3}^{2}h_{2}^{3}h_{-2}^{3}-6\,{h_{{2}}
}^{2}h_{3}^{3}h_{-2}^{2}h_{-3}^{2}-6\,h_{2}^{8}{h_{{-3
}}}^{3}h_{-2}^{2}+6\,h_{2}^{3}h_{-2}^{3}h_{{-3}}h_{3}^{2}
\n \\
&\quad
+17\,h_{2}^{6}h_{-3}^{4}h_{{3}}-6\,h_{-2}^{5}h_{2}^
{5}h_{{-3}}h_{3}^{2}+6\,h_{-3}^{3}h_{2}^{7}h_{{-2}} ) p_{{1}}=0
, \n \\
Eq_{6,4}&:
  ( -12\,h_{-3}^{2}h_{3}^{4}h_{2}^{2}h_{-2}^{5}+
108\,h_{-3}^{4}h_{2}^{5}h_{-2}^{2}h_{3}^{2}  ) r_{{1,-2}}
\n \\
&\quad
+ ( -12\,h_{2}^{10}h_{-3}^{4}h_{-2}^{3}
+42\,h_{2}^{9}h_{-3}^{4}h_{-2}^{2}+34\,h_{2}^{8}h_{-3}^{5}h_{{-2}}h_{{3}}
-30\,h_{2}^{8}h_{{-2}}h_{-3}^{4}
\n \\
&\quad
-12\,h_{2}^{7}h_{-2}^{6}h_{-3}^{2}h_{3}^{2}
-85\,h_{2}^{7}h_{{3}}h_{-3}^{5}
+24\,h_{-3}^{2}h_{2}^{6}h_{-2}^{5}h_{3}^{2}
\n \\
&\quad
+34\,h_{2}^{5}h_{-3}^{3}h_{3}^{3}h_{-2}^{4}
-12\,h_{2}^{5}h_{-3}^{2}h_{3}^{2}h_{-2}^{4}
+18\,h_{3}^{5}h_{-2}^{7}h_{{-3}}h_{2}^{2}
\n \\
&\quad
-34\,h_{2}^{4}h_{-3}^{3}h_{3}^{3}h_{-2}^{3}
-18\,h_{-2}^{8}h_{3}^{4}h_{2}^{3}
-24\,h_{2}^{3}h_{3}^{4}h_{-2}^{2}h_{-3}^{4}
-27\,h_{3}^{6}h_{-2}^{5}h_{-3}^{2}
\n \\
&\quad
+18\,h_{3}^{4}h_{-2}^{7}h_{2}^{2}
+9\,h_{{-3}}h_{3}^{5}h_{-2}^{6}h_{{2}} ) p_{{1}}=0
,     \n \\
Eq_7&:
  ( -1458\,h_{2}^{9}h_{-2}^{5}h_{-3}^{4}h_{3}^{3}+
162\,h_{2}^{6}h_{-2}^{8}h_{3}^{5}h_{-3}^{2}  ) p_{{3}}
\n \\
&\quad
+  ( 108\,h_{3}^{4}h_{-2}^{10}h_{2}^{8}h_{{-3}}+
894\,h_{2}^{9}h_{-2}^{5}h_{-3}^{4}h_{3}^{3}+294\,{h_{{
2}}}^{6}h_{-2}^{8}h_{3}^{5}h_{-3}^{2}
\n \\
&\quad
-762\,h_{2}^{8}{h
_{{3}}}^{3}h_{-3}^{4}h_{-2}^{4}+452\,h_{2}^{8}h_{3}^{4
}h_{-3}^{5}h_{-2}^{4}-72\,h_{2}^{9}h_{-2}^{5}{h_{{-3}}
}^{3}h_{3}^{2}+216\,h_{2}^{10}h_{-2}^{6}h_{-3}^{3}{h_{
{3}}}^{2}
\n \\
&\quad
-192\,h_{2}^{7}h_{-2}^{9}h_{3}^{5}h_{-3}^{2}+
306\,h_{2}^{11}h_{-3}^{7}h_{3}^{2}h_{{-2}}-216\,h_{2}^
{11}h_{-3}^{3}h_{-2}^{7}h_{3}^{2}
\n \\
&\quad
-216\,h_{-2}^{9}h_{3}^{4}h_{2}^{7}h_{{-3}}
+72\,h_{2}^{12}h_{-3}^{3}{h_{{-2}
}}^{8}h_{3}^{2}+114\,h_{2}^{5}h_{-2}^{7}h_{3}^{5}{h_{{
-3}}}^{2}
\n \\
&\quad
-348\,h_{2}^{10}h_{-2}^{6}h_{3}^{3}h_{-3}^{4}
+1092\,h_{2}^{12}h_{-3}^{6}h_{{3}}h_{-2}^{2}-312\,{h_{{2}}
}^{13}h_{-3}^{6}h_{{3}}h_{-2}^{3}
\n \\
&\quad
+72\,h_{3}^{8}h_{-2}^
{9}h_{{-3}}h_{{2}}-216\,h_{2}^{6}h_{3}^{5}h_{-3}^{6}{h_{{-
2}}}^{2}+108\,h_{3}^{4}h_{-2}^{8}h_{2}^{6}h_{{-3}}-1130\,{
h_{{-3}}}^{5}h_{2}^{7}h_{3}^{4}h_{-2}^{3}
\n \\
&\quad
-317\,h_{2}^{
4}h_{3}^{6}h_{-3}^{3}h_{-2}^{6}+144\,h_{-2}^{10}{h_{{3
}}}^{8}h_{{-3}}h_{2}^{2}-780\,h_{2}^{11}h_{-3}^{6}h_{{3}}h
_{{-2}}
\n \\
&\quad
+374\,h_{2}^{5}h_{-2}^{7}h_{3}^{6}h_{-3}^{3}-
435\,h_{3}^{7}h_{-3}^{4}h_{-2}^{5}h_{2}^{3}+144\,{h_{{
-2}}}^{10}h_{3}^{7}h_{2}^{2}
\n \\
&\quad
+432\,h_{2}^{13}h_{-3}^{5}
h_{-2}^{3}-765\,h_{2}^{10}h_{-3}^{7}h_{3}^{2}-144\,{h_
{{3}}}^{7}h_{-2}^{11}h_{2}^{3}
\n \\
&\quad
+72\,h_{-3}^{5}h_{-2}^{5
}h_{2}^{15}-216\,h_{3}^{9}h_{-2}^{8}h_{-3}^{2}-324\,{h
_{{-3}}}^{5}h_{2}^{14}h_{-2}^{4}-180\,h_{2}^{12}h_{-3}
^{5}h_{-2}^{2} ) p_{{1}}=0. \n
\end{align}

$\;$

\vsp{1cm}

$\;$

\noi
A-2. \textbf{Section \ref{sec44}}
\begin{align}
Eq_{6,1}&:
  ( 36\,h_{2}^{5}h_{-3}^{2}+36\,h_{2}^{2}h_{3}^{2}h_{-2}^{3}  ) r_{{2,1}}
\n \\
& \quad
+ ( -10\,h_{-2}^{3}h_{3}^{3}+6\,h_{2}^{3}h_{{3}}h_{{-2
}}h_{{-4}}+6\,h_{{3}}h_{{4}}h_{{2}}h_{-2}^{3}-6\,h_{{3}}h_{{4}}{h_
{{2}}}^{2}h_{-2}^{4}-3\,h_{{2}}h_{-2}^{2}h_{4}^{2}h_{{-3}}
\n \\
& \quad
-6\,h_{2}^{4}h_{{3}}h_{-2}^{2}h_{{-4}}+10\,h_{{3}}h_{2}^{4
}h_{-3}^{2}h_{{-2}}+4\,h_{3}^{2}h_{-2}^{2}h_{{4}}h_{{-3}}+
4\,h_{2}^{2}h_{{-3}}h_{3}^{2}h_{{-4}}-3\,h_{2}^{3}h_{{-3}}
h_{{-4}}h_{{4}}
\n \\
& \quad
-10\,h_{{3}}h_{2}^{3}h_{-3}^{2}
+10\,h_{{2}}h_{-2}^{4}h_{3}^{3}  ) p_{{1}}=0,
\n \\
Eq_{6,2}&:
 ( -8\,h_{2}^{3}h_{{3}}h_{{-2}}h_{-3}^{2}
-8\,h_{3}^{3}h_{-2}^{4} ) r_{{1,-2}}
 + ( 12\,h_{2}^{5}h_{-3}^{2}h_{{-2}}+12
\,h_{2}^{2}h_{-2}^{4}h_{3}^{2} ) r_{{2,-1}}
\n \\
&  \quad
+ ( h_{2}^{3}h_{{-3}}h_{-2}^{2}h_{{4}}+2\,h_{2}^{3}h_{
{-3}}h_{3}^{2}h_{-2}^{3}-h_{2}^{6}h_{{-2}}h_{{-4}}h_{{-3}}
-h_{2}^{4}h_{{-3}}h_{-2}^{3}h_{{4}}+h_{2}^{5}h_{{-4}}h_{{-
3}}
\n \\
&\quad
+2\,h_{2}^{6}h_{-3}^{3}  ) p_{{1}}=0
\n \\
Eq_{6,3}&:
 ( 36\,h_{2}^{3}h_{-2}^{4}h_{{-3}}h_{3}^{2}+36
\,h_{2}^{6}h_{{-2}}h_{-3}^{3} ) r_{{2,-1}}
\n \\
&\quad
 +( 6\,h_{-3}^{4}h_{2}^{7}-3\,h_{{2}}h_{-2}^{4}{h_{{4}
}}^{2}h_{{3}}h_{{-3}}+3\,h_{-2}^{5}h_{{4}}h_{{2}}h_{3}^{2}-3\,
h_{-2}^{3}h_{{4}}h_{2}^{5}h_{-3}^{2}+4\,h_{-2}^{4}h_{{
4}}h_{3}^{3}h_{{-3}}
\n \\
&\quad
+3\,h_{-2}^{2}h_{{4}}h_{2}^{4}{h_{{-3}
}}^{2}+10\,h_{2}^{4}h_{3}^{2}h_{-3}^{2}h_{-2}^{3}-3\,h
_{{-4}}h_{2}^{4}h_{3}^{2}h_{-2}^{4}-10\,h_{3}^{4}{h_{{
-2}}}^{5}+4\,h_{2}^{3}h_{{4}}h_{-3}^{3}h_{{3}}h_{{-2}}
\n \\
&\quad
+4\,h_{3}^{4}h_{-2}^{6}h_{{2}}-10\,h_{-2}^{2}h_{2}^{3}{h_{{-3}
}}^{2}h_{3}^{2}-3\,h_{2}^{3}h_{{-3}}h_{-2}^{2}h_{{3}}h_{{-
4}}h_{{4}}-3\,h_{{-4}}h_{2}^{7}h_{-3}^{2}h_{{-2}}
\n \\
&\quad
+3\,h_{2}^{6}h_{{-4}}h_{-3}^{2}
+3\,h_{{-4}}h_{2}^{3}h_{3}^{2}h_{-2}^{3}
-3\,h_{-2}^{6}h_{{4}}h_{2}^{2}h_{3}^{2}  ) p_{{1}}=0,
\n \\
Eq_{7}&:
 ( 216\,h_{{-2}}h_{2}^{8}h_{-3}^{4}+432\,h_{3}^{2
}h_{-2}^{4}h_{2}^{5}h_{-3}^{2}+216\,h_{3}^{4}{h_{{-2}}
}^{7}h_{2}^{2}  ) p_{{3}}
\n \\
&\quad
+ ( 20\,h_{3}^{4}h_{-2}^{6}h_{{2}}+40\,h_{{3}}h_{2}^{7
}h_{-3}^{5}-36\,h_{2}^{4}h_{{-4}}h_{{3}}h_{-2}^{3}h_{{4}}h
_{{-3}}+18\,h_{2}^{2}h_{{3}}h_{-2}^{5}h_{4}^{2}h_{{-3}}
\n \\
&\quad
+18\,h_{2}^{4}h_{-2}^{2}h_{4}^{2}h_{-3}^{2}h_{{-4}}-12\,{
h_{{2}}}^{3}h_{-4}^{2}h_{{4}}h_{3}^{2}h_{-2}^{3}-12\,h_{{2
}}h_{-2}^{5}h_{4}^{2}h_{{-4}}h_{3}^{2}-60\,h_{-3}^{3}{
h_{{2}}}^{7}h_{{3}}h_{{-2}}h_{{-4}}
\n \\
&\quad
-68\,h_{3}^{3}h_{-2}^{3}h_{{-4}}h_{2}^{3}h_{{-3}}
+12\,h_{3}^{2}h_{-2}^{2}h_{{4}}{h_{{
-3}}}^{2}h_{2}^{3}h_{{-4}}+32\,h_{3}^{2}h_{{-2}}h_{2}^{3}h
_{{4}}h_{-3}^{4}
\n \\
&\quad
+8\,h_{3}^{2}h_{-2}^{4}h_{4}^{2}{h_{{-
3}}}^{2}h_{{2}}-48\,h_{3}^{4}h_{-2}^{2}h_{2}^{2}h_{-3}
^{2}h_{{-4}}+44\,h_{2}^{4}h_{3}^{3}h_{-2}^{4}h_{{-4}}h_{{-
3}}+54\,h_{{3}}h_{-2}^{2}h_{2}^{4}h_{-3}^{3}h_{{4}}
\n \\
&\quad
-66\,h_{-3}^{3}h_{2}^{5}h_{-2}^{3}h_{{3}}h_{{4}}+36\,h_{2}^{5
}h_{{3}}h_{-2}^{4}h_{{4}}h_{{-4}}h_{{-3}}-18\,h_{4}^{2}h_{{3}}
h_{-2}^{6}h_{2}^{3}h_{{-3}}-16\,h_{3}^{5}h_{-2}^{6}h_{
{2}}h_{{-3}}
\n \\
&\quad
+24\,h_{3}^{3}h_{-2}^{3}h_{2}^{4}h_{-3}^{3
}+38\,h_{2}^{2}h_{3}^{3}h_{-2}^{6}h_{{4}}h_{{-3}}+24\,h_{{
3}}h_{{-2}}h_{2}^{6}h_{{-5}}h_{-3}^{2}+54\,h_{2}^{7}h_{{3}
}h_{-2}^{2}h_{-4}^{2}h_{{-3}}
\n \\
&\quad
+16\,h_{3}^{4}h_{-2}^{5}h
_{{4}}h_{{-4}}-20\,h_{2}^{5}h_{-3}^{2}h_{3}^{2}h_{-4}^
{2}-24\,h_{{-2}}h_{2}^{4}h_{4}^{2}h_{-3}^{4}-16\,h_{3}
^{4}h_{-2}^{4}h_{{4}}h_{-3}^{2}
\n \\
&\quad
-12\,h_{3}^{2}h_{-2}^{3
}h_{2}^{3}h_{{-3}}h_{{4}}h_{{-5}}-50\,h_{3}^{3}h_{-2}^{5}h
_{{4}}h_{{-3}}h_{{2}}-24\,h_{{3}}h_{-2}^{2}h_{2}^{7}h_{{-5}}{h
_{{-3}}}^{2}+24\,h_{-2}^{3}h_{2}^{9}h_{{-4}}h_{-3}^{2}
\n \\
&\quad
+16\,h_{2}^{2}h_{3}^{4}h_{-4}^{2}h_{-2}^{3}+15\,h_{2}
^{6}h_{-3}^{2}h_{-4}^{2}h_{{4}}+40\,h_{3}^{3}h_{-2}^{2
}h_{2}^{3}h_{-3}^{3}+16\,h_{3}^{2}h_{2}^{5}h_{-3}^
{3}h_{{-5}}
\n \\
&\quad
+36\,h_{{3}}h_{2}^{6}h_{-3}^{3}h_{{-4}}+24\,{h_{{3}
}}^{3}h_{-2}^{4}h_{2}^{3}h_{{-5}}-48\,h_{-2}^{2}h_{2}^
{8}h_{{-4}}h_{-3}^{2}-24\,h_{3}^{3}h_{-2}^{5}h_{2}^{4}
h_{{-5}}
\n \\
&\quad
-54\,h_{{3}}h_{{-2}}h_{2}^{6}h_{-4}^{2}h_{{-3}}+48\,{h
_{{3}}}^{2}h_{-2}^{7}h_{2}^{3}h_{{4}}-48\,h_{-2}^{4}{h_{{2
}}}^{6}h_{{4}}h_{-3}^{2}-24\,h_{3}^{2}h_{-2}^{6}h_{2}^
{6}h_{{-4}}
\n \\
&\quad
+48\,h_{3}^{2}h_{-2}^{5}h_{2}^{5}h_{{-4}}-24\,{
h_{{3}}}^{2}h_{-2}^{8}h_{2}^{4}h_{{4}}+52\,h_{3}^{2}{h_{{-
2}}}^{5}h_{2}^{6}h_{-3}^{2}-32\,h_{3}^{2}h_{-2}^{4}{h_
{{2}}}^{5}h_{-3}^{2}
\n \\
&\quad
+24\,h_{-2}^{5}h_{2}^{7}h_{{4}}{h_{{-3
}}}^{2}-12\,h_{2}^{6}h_{-3}^{3}h_{{4}}h_{{-5}}+16\,h_{3}^{
4}h_{-2}^{3}h_{2}^{2}h_{{-3}}h_{{-5}}-24\,h_{-2}^{6}h_{{4}
}h_{2}^{2}h_{3}^{2}
\n \\
&\quad
+3\,h_{2}^{2}h_{-2}^{4}h_{4}^{3
}h_{-3}^{2}+24\,h_{{-4}}h_{2}^{7}h_{-3}^{2}h_{{-2}}+24\,{h
_{{-2}}}^{3}h_{{4}}h_{2}^{5}h_{-3}^{2}-20\,h_{2}^{4}{h_{{3
}}}^{2}h_{-3}^{2}h_{-2}^{3}
\n \\
&\quad
-24\,h_{{-4}}h_{2}^{4}h_{3}
^{2}h_{-2}^{4}-40\,h_{-3}^{4}h_{2}^{7}-16\,h_{2}^{9}{h
_{{-3}}}^{4}h_{-2}^{2}+40\,h_{3}^{5}h_{-2}^{5}h_{{-3}}+56
\,h_{{-2}}h_{2}^{8}h_{-3}^{4}
\n \\
&\quad
+68\,h_{3}^{4}h_{-2}^{8}{
h_{{2}}}^{3}-88\,h_{3}^{4}h_{-2}^{7}h_{2}^{2}  ) p_{{1}}=0.
\n
\end{align}

$\;$

\vsp{1.5cm}

$\;$

\noi
\texttt{Step 8} in  Section \ref{sec442},  formulae for $\, h_{\pm 5}$:
\begin{align}
h_5&=
- \, \frac{1}{12 \,h_{2}^{4}h_{-3}^{3}h_{{-2}}
( h_{2}^{3}h_{-3}^{2}+h_{-2}^{3}h_{3}^{2} ) }
(
 36\,h_{-3}^{3}h_{2}^{5}h_{{-4}}h_{-2}^{2}{h
_{{3}}}^{3}+12\,h_{2}^{4}h_{-3}^{2}h_{3}^{2}h_{-2}^{5}
h_{{4}}
\n \\
&\quad
-10\,h_{3}^{6}h_{-2}^{8}-30\,h_{2}^{4}h_{-3}^{3
}h_{-2}^{4}h_{4}^{2}h_{{3}}-12\,h_{2}^{7}h_{3}^{2}{h_{
{-2}}}^{4}h_{{-4}}h_{-3}^{2}
\n \\
&\quad
-30\,h_{2}^{10}h_{-3}^{6}-12\,
h_{2}^{5}h_{-2}^{6}h_{3}^{2}h_{{4}}h_{-3}^{2}+20\,{h_{
{3}}}^{4}h_{-2}^{6}h_{2}^{4}h_{-3}^{2}-14\,h_{2}^{7}{h
_{{3}}}^{2}h_{-2}^{3}h_{-3}^{4}
\n \\
&\quad
-40\,h_{2}^{6}h_{-3}^{4
}h_{3}^{2}h_{-2}^{2}-15\,h_{2}^{7}h_{-3}^{4}h_{-2}
^{2}h_{{4}}+4\,h_{3}^{5}h_{-2}^{7}h_{{4}}h_{{-3}}+15\,{h_{{2}}
}^{10}h_{{-2}}h_{{-4}}h_{-3}^{4}
\n \\
&\quad
+3\,h_{2}^{3}h_{3}^{4}{h_{
{-2}}}^{6}h_{{-4}}+3\,h_{{2}}h_{3}^{4}h_{-2}^{8}h_{{4}}-50\,{h
_{{2}}}^{3}h_{-3}^{2}h_{3}^{4}h_{-2}^{5}+15\,h_{2}^{8}
h_{-2}^{3}h_{{4}}h_{-3}^{4}
\n \\
&\quad
-3\,h_{2}^{2}h_{3}^{4}{h_{{
-2}}}^{9}h_{{4}}-3\,h_{2}^{4}h_{3}^{4}h_{-2}^{7}h_{{-4}}-
15\,h_{2}^{9}h_{{-4}}h_{-3}^{4}+4\,h_{{2}}h_{3}^{6}{h_{{-2
}}}^{9}
\n \\
&\quad
+12\,h_{2}^{6}h_{-3}^{2}h_{3}^{2}h_{-2}^{3}h_{{
-4}}+20\,h_{{4}}h_{-3}^{3}h_{2}^{3}h_{3}^{3}h_{-2}^{4}
-20\,h_{{4}}h_{-3}^{5}h_{2}^{6}h_{{3}}h_{{-2}}
\n \\
&\quad
-3\,h_{{2}}h_{{-
3}}h_{-2}^{7}h_{4}^{2}h_{3}^{3}-30\,h_{2}^{6}{h_{{-3}}
}^{3}h_{-2}^{2}h_{{-4}}h_{{4}}h_{{3}}+3\,h_{2}^{5}h_{-2}^{
3}h_{3}^{3}h_{-4}^{2}h_{{-3}}
),
\n \\
h_{-5}&=
 \, \frac{1}{12 h_{2}^{2}h_{{3}}h_{{-2}} ( h_{2}^
{3}h_{-3}^{2}+h_{-2}^{3}h_{3}^{2}  ) }
(
-12\,h_{2}^{2}h_{-2}^{7}h_{{4}}h_{3}^{2}+30
\,h_{{-3}}h_{2}^{3}h_{{3}}h_{-2}^{3}h_{{-4}}h_{{4}}
\n \\
&\quad
+12\,h_{2}^{3}h_{-2}^{4}h_{{-4}}h_{3}^{2}+50\,h_{-3}^{2}h_{2}
^{3}h_{3}^{2}h_{-2}^{3}-4\,h_{3}^{3}h_{-2}^{5}h_{{4}}h
_{{-3}}-40\,h_{{3}}h_{-2}^{2}h_{-3}^{3}h_{{4}}h_{2}^{3}
\n \\
&\quad
+10\,h_{3}^{4}h_{-2}^{6}-12\,h_{2}^{4}h_{-2}^{5}h_{{-4}}{
h_{{3}}}^{2}-12\,h_{2}^{4}h_{-2}^{3}h_{{4}}h_{-3}^{2}+3\,h
_{{2}}h_{-2}^{5}h_{4}^{2}h_{{3}}h_{{-3}}+10\,h_{-3}^{2}{h_
{{2}}}^{4}h_{-2}^{4}h_{3}^{2}
\n \\
&\quad
+27\,h_{{-3}}h_{2}^{5}h_{{3}}
h_{{-2}}h_{-4}^{2}+40\,h_{2}^{6}h_{-3}^{4}+12\,h_{2}^{
5}h_{-2}^{4}h_{{4}}h_{-3}^{2}-12\,h_{2}^{6}h_{{-2}}h_{{-4}
}h_{-3}^{2}-16\,h_{-3}^{3}h_{2}^{5}h_{{3}}h_{{-4}}
\n \\
&\quad
+12\,h_{
{2}}h_{-2}^{6}h_{{4}}h_{3}^{2}+12\,h_{2}^{7}h_{-2}^{2}
h_{{-4}}h_{-3}^{2}-16\,h_{-3}^{4}h_{2}^{7}h_{{-2}}+26\,{h_
{{3}}}^{4}h_{-2}^{7}h_{{2}}+20\,h_{{-3}}h_{2}^{2}h_{3}^{3}
h_{{-4}}h_{-2}^{3}
). \n
\end{align}
\noi
Note that $\, h_{-5} \ne \overline{h_5}$ in this  formula.
$\, h_{-5} -\overline{h_5}=0$ gives an integrability equation
which is quadratic in $\, h_{4}, \, h_{-4}$.


\end{document}